%% file: agt.tex
\documentclass{amsart}     
\usepackage{amssymb,amsthm,color}
\usepackage{epsfig}
\usepackage[all]{xy}

\title{Actions of certain arithmetic groups on Gromov hyperbolic spaces}

\author[J. F. Manning]{Jason Fox Manning}
\address{University at Buffalo, SUNY}
\email{j399m@buffalo.edu}
\urladdr{http://www.math.buffalo.edu/~j399m}

\newtheorem{theorem}{Theorem}[section]
\newtheorem{lemma}[theorem]{Lemma}
\newtheorem{corollary}[theorem]{Corollary}
\newtheorem{proposition}[theorem]{Proposition}
\newtheorem{question}[theorem]{Question}
\newtheorem{conjecture}[theorem]{Conjecture}

\newtheorem{case}{Case}
\newtheorem{subcase}{Case}[case]

\newtheorem{claim}[theorem]{Claim}
\newtheorem{subclaim}[theorem]{Subclaim}

\theoremstyle{definition}
\newtheorem{remark}[theorem]{Remark}
\newtheorem{definition}[theorem]{Definition}

\newtheorem{observation}[theorem]{Observation}

 \newcommand{\cd}[1]{\begin{equation*}{\xymatrix{#1}}\end{equation*}}
 \newcommand{\cdlabel}[2]{\begin{equation}\label{#1}{\xymatrix{#2}}\end{equation}}
\newcommand{\FA}{\ensuremath{(\mathrm{FA})}}
\newcommand{\QFA}{\ensuremath{(\mathrm{QFA})}}

\newcommand{\Mtwo}[4]{\ensuremath{\begin{pmatrix} {#1} & {#2} \\ {#3} &
    {#4}\end{pmatrix} }}
\newcommand{\textMtwo}[4]{\ensuremath{\left(\begin{smallmatrix} {#1} & {#2} \\ {#3} &
    {#4}\end{smallmatrix}\right) }}
\newcommand{\mc}[1]{\mathcal{#1}}
\newcommand{\uu}[1]{\ensuremath{\underline{#1}}}
\newcommand{\gp}[3]{\ensuremath{({#1}\cdot{#2})_{#3}}}

\def\co{\colon\thinspace}
\def\R{\mathbb{R}}
\def\Z{\mathbb{Z}}
\def\C{\mathbb{C}}
\def\H{\mathbb{H}}
\def\N{\mathbb{N}}
\def\Q{\mathbb{Q}}

\begin{document}

\begin{abstract}    
We study the variety of actions of a fixed (Chevalley) group on
arbitrary geodesic, Gromov hyperbolic spaces.  In high rank we obtain
a complete classification. In rank one, we obtain some partial results
and give a conjectural picture.
\end{abstract}

\maketitle

\tableofcontents

\section{Introduction}
Given a group $G$ one may ask the question: 
\begin{question}
In what ways can $G$ act non-trivially on a Gromov hyperbolic metric space?
\end{question}
Many interesting groups can be fruitfully studied via some natural action on
a Gromov hyperbolic space.  Examples include the action of an
amalgam or HNN extension on its Bass-Serre tree, the action of a
Kleinian group on $\H^n$, and the action of the mapping class group of
a surface on the curve complex of that same surface.  Alternatively, one can
study the space of actions of a fixed group on some (fixed or
varying) Gromov hyperbolic metric space.  For example, the $PSL(2,\C)$--character variety of a group
parameterizes the space of actions of a fixed
group on the hyperbolic space $\H^3$.  Analysis of the structure of
this variety has led to many remarkable theorems about $3$--manifolds
and their fundamental
groups (for an introduction, see \cite{shalen:handbook}).
A larger ``variety'' (in scare quotes because there
is unlikely to be any algebraic structure) would describe all
non-trivial actions on Gromov hyperbolic spaces, up to some
appropriate equivalence relation.  We give 
suggestions for how to define this equivalence relation and topologize the
variety in Section \ref{s:actions}.  Briefly, the equivalence relation
is that generated by coarsely equivariant quasi-isometric embeddings.
This equivalence is coarser than that given by quasi-conjugacy (as in
for instance \cite{msw:quasiactI}), but finer than that given by
equivariant homeomorphism of limit sets.
In this paper we concentrate on cases in which the set of equivalence
classes is particularly simple.

Certain equivalence classes of actions on hyperbolic spaces cannot be
ruled out, or 
even really analyzed using the tools of negative curvature.  These are
the \emph{actions with an invariant horoball} (see Section
\ref{ss:combhoro} for the definition, and Theorem \ref{t:horeq} for
some characterizations).  Actions with an invariant
horoball (in the sense used in this paper) are always elementary; they
include the trivial action
on a point and actions which preserve some horofunction.  
A cobounded
action on an unbounded space never has an invariant horoball.

The variety of actions of an irreducible higher rank lattice is
expected to be very simple.
\begin{conjecture}\label{c:main}
If $\Gamma$ is an irreducible lattice in a higher rank Lie group (or
in a 
nontrivial product of locally compact groups) $G$,
there there are finitely many Gromov hyperbolic $G$--spaces (up to
coarse equivalence) without
invariant horoballs.
\end{conjecture}
In the case where $G$ is a simple Lie group of rank at least $2$, the
only expected actions are those with an invariant horoball.  If $G$
has more than one direct factor, then $\Gamma$ projects densely to
each factor.  If the factors have rank one, then there will clearly be
non-elementary isometric actions of $\Gamma$ on rank one symmetric
spaces.  We discuss this case in Section \ref{section:rankone}.

For actions by lattices in nontrivial products of simple Lie
groups, the conjecture follows from rigidity results of Monod
and Monod-Shalom
\cite{monod:superrigid,monodshalom:jdg} 
if one restricts attention
to Gromov hyperbolic spaces which are also either CAT$(0)$ spaces,
proper and cocompact spaces, or bounded valence graphs.
(Cf.
\cite{gelanderkarlssonmargulis} for CAT$(0)$ spaces.)  
For an example
of an action of a lattice on a Gromov hyperbolic space which is
inequivalent to any action on a Gromov hyperbolic CAT$(0)$ space,
consider a lattice as in the Appendix to \cite{manning:qfa}, which has
Property (T), but admits a non-trivial pseudocharacter (or homogeneous
quasi-morphism).  This
pseudocharacter gives rise to a cobounded action on a space
quasi-isometric to $\R$, fixing both ends.  On the other hand, an
action on a CAT$(0)$ space preserving a point at infinity would give
(via the Busemann function) a homomorphism to $\R$, and  a cobounded such
action would give an unbounded homomorphism to $\R$.
Such a homomorphism is ruled out by Property (T).

\begin{theorem}\label{thm:cdr2} 
Suppose that $G$ is a simple Chevalley-Demazure
group scheme of rank at least $2$, and let $\mathcal{O}$ be the ring
of integers of any number field.  Then any isometric action of
$G(\mathcal{O})$ on a Gromov hyperbolic geodesic metric space has
an invariant horoball.
\end{theorem}
Some remarks:
\begin{enumerate}
\item Some closely related results are proved by Karlsson and Noskov
  \cite[Sections 8.2 and 8.3]{karlssonnoskov}; part of our strategy is
  similar to theirs, and to that of Fukunaga in \cite{fukunaga:FA}.
\item The rank $\geq 2$ assumption is necessary, as the
  action of $SL(2,\mathcal{O})$ on $\H^3$ obtained from the inclusion 
  $SL(2,\mathcal{O})\longrightarrow SL(2,\C)$  never has an
  invariant horoball.
\end{enumerate}

The rank one case is discussed in Section 
\ref{section:rankone}, where we apply a result of Carter, Keller and
Paige to show a weaker theorem for these groups.
\begin{theorem}\label{thm:cdr1}
Let $\mathcal{O}$ be the ring of integers of a number field, and
suppose that $\mathcal{O}$ has infinitely many units.  Suppose $X$ is
quasi-isometric to a tree.  Every action of $SL(2,\mathcal{O})$ on $X$
has a bounded orbit.
\end{theorem}

\subsection{Relative hyperbolicity and bounded generation}
If $G$ is a finitely generated
group, $S$ is a generating set for $G$, and
$\mc{P}=\{P_1,\ldots,P_n\}$ is a collection of
subgroups of $G$, then one can form the \emph{coned space}
$C(G,\mc{P},S)$ as follows:  Let $\Gamma(G,S)$ be the Cayley graph of
$G$ with respect to $S$. 
The coned space $C(G,\mc{P},S)$ is obtained
from $\Gamma(G,S)$ by coning each left coset of an element of $\mc{P}$
to a point.  (Here the $0$--skeleton of $\Gamma(G,S)$ is implicitly
identified with $G$.)

Recall that a graph is \emph{fine} if every edge is contained in only
finitely many circuits (i.e. embedded cycles) of any bounded length.
\begin{definition}
If $G$ is a group with finite generating set $S$, and a collection of
subgroups $\mc{P}=\{P_1,\ldots,P_n\}$, then
$G$ is \emph{weakly hyperbolic relative to} $\mc{P}$ if
$C(G,\mc{P},S)$ is $\delta$--hyperbolic for some $\delta$.  

If, moreover, $C(G,\mc{P},S)$ is \emph{fine}, then $G$ is
\emph{(strongly) hyperbolic relative to} $\mc{P}$. 
\end{definition}

A special case of weak relative hyperbolicity is bounded generation.
The following definition is easily seen to be equivalent to the
standard one:
\begin{definition}
A group $G$ is \emph{boundedly generated} by a collection of
subgroups $\mc{P}=\{P_1,\ldots,P_n\}$ if the Cayley graph of $G$ with
respect to $\cup\mc{P}$ has finite diameter.

If each $P_i$ is cyclic, generated by $g_i$, we say that $G$ is
\emph{boundedly generated by} $\{g_1,\ldots,g_n\}$.
\end{definition}

A cobounded action on an unbounded Gromov hyperbolic space does not
have an invariant horoball.
It is thus a corollary of Theorem \ref{thm:cdr2} that these higher rank
$G(\mathcal{O})$ are not strongly relatively hyperbolic
with respect to any system of proper subgroups.  
\begin{corollary}\label{corollary:notrelhyp}
If $L=G(\mc{O})$ is as in Theorem \ref{thm:cdr2}, and $L$ is weakly
hyperbolic relative to a system of subgroups
$\mc{P}$, then the coned space $C(L,\mc{P},S)$ has finite
diameter for any finite generating set $S$.

In particular, $L$ is not strongly hyperbolic relative to any system
of proper subgroups.
\end{corollary}
By possibly altering $\mc{P}$ to add some cyclic subgroups, the first
part of the corollary can be restated:  If $L$ is weakly 
hyperbolic relative to a system of subgroups $\mc{P}$ which generate
$L$, then $L$ is boundedly generated by $\mc{P}$.
\begin{remark}
Corollary \ref{corollary:notrelhyp} also can be easily deduced from
\cite{karlssonnoskov} in some special cases, including $L=SL(n,\Z)$
for $n>2$.  The second part of Corollary \ref{corollary:notrelhyp}
(about strong relative hyperbolicity) can be deduced from known
theorems in a number of ways, perhaps most straightforwardly by
combining a theorem of Fujiwara \cite{fujiwara:gromovhyperbolic} with
one of Burger and Monod \cite{burgermonod:gafa}.
\end{remark}

\subsection{Outline}
In Section \ref{s:preliminaries} we recall some definitions and basic
results, first from the theory of Gromov hyperbolic spaces, and second
from Chevalley groups over number rings.  In Section \ref{s:actions}
an equivalence relation amongst hyperbolic $G$--spaces is proposed, and
generalized combinatorial horoballs are introduced.  In Section
\ref{s:elementary}, we improve on the statement and proof of
Proposition 3.9 of \cite{manning:qfa}, and use the improved version to
characterize hyperbolic $G$--spaces with invariant horoballs.
In Section \ref{s:bigrank} we prove Theorem \ref{thm:cdr2}, and in
Section \ref{section:rankone} we prove Theorem \ref{thm:cdr1}.

\subsection{Acknowledgments}
The author thanks Benson Farb and Hee Oh for useful conversations, and
Nicolas Monod for helpful comments on an earlier version of this paper.
This work was partly supported by an NSF Postdoctoral Research
Fellowship (DMS-0301954). 

\section{Preliminaries}\label{s:preliminaries}
\subsection{Coarse geometry}
\begin{definition}
If $X$ and $Y$ are metric spaces, $K\geq 1$ and $C\geq 0$,
a \emph{$(K,C)$--quasi-isometric embedding} of $X$ into $Y$ is a function
$q\co X\to Y$ so that
 For all $x_1$, $x_2\in X$
\[\frac{1}{K}d(x_1,x_2)-C\leq d(q(x_1),q(x_2))\leq Kd(x_1,x_2)+C\]

If in addition the map $q$ is \emph{$C$--coarsely onto} -- \emph{i.e.},
every $y\in Y$ is
distance at most $C$ from some point in $q(X)$ -- then $q$ is called a
\emph{$(K,C)$--quasi-isometry}.
The two metric spaces $X$ and $Y$ are then
said to be \emph{quasi-isometric} to one
another.  This is a symmetric condition.
\end{definition}
\begin{definition}
A \emph{$(K,C)$--quasi-geodesic} in $X$
is a $(K,C)$--quasi-isometric embedding $\gamma\co\R\to X$.
We will occasionally abuse
notation by referring to the image of $\gamma$ as a quasi-geodesic.
\end{definition}
\subsection{Gromov hyperbolic spaces}
For more details on Gromov hyperbolic metric spaces, see
\cite{bridhaef:book} and \cite{gromov:wordhyperbolic}.
A number of equivalent definitions of Gromov hyperbolicity are known.
For geodesic spaces, we will use the one based on the existence of
\emph{comparison tripods}. 
Given a geodesic
triangle $\Delta(x,y,z)$ in
any metric space, there is a unique \emph{comparison tripod}, $T_\Delta$,
a metric tree so
that the distances between the three extremal points of the tree,
$\overline{x}$, $\overline{y}$ and $\overline{z}$ , are the same as
the distances between $x$, $y$ and $z$.
There is an obvious map $\pi\co \Delta(x,y,z)\to T_\Delta$ which takes $x$ to
$\overline{x}$, $y$ to $\overline{y}$ and $z$ to $\overline{z}$, and
which is an isometry on each side of $\Delta(x,y,z)$.
\begin{definition}\label{d:geodgromov}
A geodesic space $X$ is \emph{$\delta$--hyperbolic} if for any geodesic triangle
$\Delta(x,y,z)$ and any point $p$ in the comparison tripod $T_\Delta$,
 the diameter of $\pi^{-1}(p)$ is less than $\delta$.  If $\delta$ is
 unimportant we may simply say that $X$ is \emph{Gromov hyperbolic}.
\end{definition}
Gromov hyperbolicity (of geodesic spaces) is a quasi-isometry invariant.

The following proposition about stability of quasi-geodesics in Gromov
hyperbolic spaces
is well-known (see, e.g. \cite[III.H.1.7]{bridhaef:book}).
\begin{proposition}\label{p:quasistable}
Let $K\geq 1$, $C\geq 0$, $\delta\geq 0$.  Then there is some
$B=B(K,C,\delta)$, so that whenever $\gamma$ and $\gamma'$ are two
$(K,C)$--quasi-geodesics with the same endpoints in a
$\delta$--hyperbolic geodesic metric space $X$, then the Hausdorff
distance between $\gamma$ and $\gamma'$ is at most $B$.
\end{proposition}

We will occasionally need to deal with spaces which are not geodesic.
If $X$ is $\delta$--hyperbolic in the sense of Definition
\ref{d:geodgromov}, then it
satisfies the four point condition:  For all $p_1$, $p_2$, $p_3$,
$p_4\in X$, 
\begin{equation}\label{fourpoint}
d(p_1,p_4)+d(p_2,p_3)\leq\max\{d(p_1,p_2)+d(p_3,p_4)\mbox{, }d(p_1,p_3)+d(p_2,p_4)+2\delta\}.
\end{equation}
Conversely, if a geodesic space satisfies \eqref{fourpoint}, then it
is $6\delta$--hyperbolic in the sense of Definition \ref{d:geodgromov}.
(For both these facts, and the below definition, see
\cite[III.H.1.22]{bridhaef:book} or \cite{gromov:wordhyperbolic}.)
\begin{definition}\label{d:fourpoint}
A (not necessarily geodesic) metric space $X$ is \emph{$(\delta)$--hyperbolic}
if it satisfies the condition \eqref{fourpoint} above.  If $\delta$ is
unimportant, we simply say that $X$ is \emph{Gromov hyperbolic}.
\end{definition}

In order to describe the boundary of a Gromov hyperbolic space, we
introduce the \emph{Gromov product} notation.
\begin{definition}
If $x,y$ and $z$ are points in a metric space with a metric
$d(\cdot,\cdot)$, then
\[ \gp{x}{y}{z}:= \frac{1}{2}(d(x,z)+d(y,z)-d(x,y)).\]
\end{definition}
We should remark that in the context of Definition \ref{d:geodgromov},
$\gp{x}{y}{z}$ is the distance from
$\overline{z}$ to the central vertex of the comparison tripod for a
geodesic triangle with vertices $x$, $y$ and $z$.
\begin{definition}
Let $p\in X$ where $X$ is a Gromov hyperbolic space.
A sequence of points $\{x_i\}$ in a Gromov hyperbolic space
\emph{tends to infinity} if
$\lim_{i,j\to\infty}\gp{x_i}{x_j}{p}=\infty$.
Two such sequences are \emph{equivalent}, written $\{x_i\}\sim\{y_i\}$, if
$\lim_{i,j\to\infty}\gp{x_i}{y_j}{p}=\infty$.
The boundary $\partial X$ is the set of equivalence classes of
sequences which tend to infinity.
\end{definition}
The Gromov product extends (by taking a lim sup) to sequences which
tend to infinity, and this allows convergence in $\partial X$ to be
defined, giving $\partial X$ a natural topology.

\subsection{Isometries of Gromov hyperbolic spaces}
Isometries of geodesic hyperbolic spaces can be classified into three
types.
\begin{definition}
Let $f\co X\to X$ be an isometry.  If $x\in X$, we let
$O_x=\{f^n(x)\mid n\in \Z\}$.
We say that $f$ is \emph{elliptic} if $O_x$ is bounded.  We say that
$f$ is \emph{hyperbolic} if
$n\mapsto f^{n}(x)$ is a quasi-isometric embedding of $\Z$ into $X$.
We say that $f$ is \emph{parabolic} if $O_x$ has a unique limit point
in $\partial G$.
\end{definition}
The following was observed by Gromov \cite[8.1.B]{gromov:wordhyperbolic}.
Although it is often stated with an extra hypothesis of
properness, this hypothesis is unnecessary (See, for example the
proof in \cite[Chapitre 9]{cdp}, where the extra hypothesis is given,
but not used).
\begin{proposition}\label{p:isometryclass}
Every isometry of a geodesic Gromov hyperbolic space is elliptic,
parabolic, or hyperbolic.
\end{proposition}

\begin{lemma}\label{l:commute}
Suppose that $G$ acts on the geodesic Gromov hyperbolic space $X$, and that
$p\in G$ acts parabolically, fixing $e\in \partial X$.  If $g\in G$
commutes with $p$, then $g$ also fixes $e$.
\end{lemma}
\begin{proof}
Let $e'=g(e)$.  We have $e'=gp(e)=pg(e)=p(e')$, so $p$ fixes $e'$.
Since $p$ is parabolic, it fixes a unique point in $\partial X$, and
so $e'=e$.
\end{proof}

\begin{definition}
Let $G$ be a finitely generated group, and let $g\in G$.  If $n\mapsto
g^n$ is a quasi-isometric embedding, we say that $g$ is
\emph{undistorted}.  Otherwise, $g$ is \emph{distorted}.
\end{definition}
The proof of the following observation is left to the reader.
\begin{lemma}\label{lemma:distortion}
Let the finitely generated group $G$ act by isometries on a Gromov
hyperbolic space $X$.  If $g\in 
G$ acts hyperbolically on $X$, then $g$ is undistorted.
\end{lemma}

\subsection{Chevalley groups}\label{defs:chevalley}
In this section we recall the definition of a Chevalley group over a
commutative ring.
 (All rings are assumed to have $1 \neq 0$.)
The simplest example of a Chevalley group is
$SL(n,\Z)$.  If one thinks of $SL(n,\Z)$ as being the ``$\Z$--points
of $SL(n,\C)$,'' then the Chevalley-Demazure group scheme
identifies what the ``$R$--points of $G$'' are, where
$R$ is now allowed to be an arbitrary ring, and $G$ an arbitrary
complex semisimple Lie group.  It turns out that this idea is not entirely
well-formed, until one fixes an embedding of $G$ into $GL(n,\C)$ for
some $n$.  The following exposition is largely 
adapted from \cite{abe:chevalley} and \cite{tavgen:bg}.

Let $\rho\co G\to GL(V)$ be a representation of a connected complex semisimple
Lie group into the general linear group of a complex vector space $V$
of dimension $n$.  We will assume that
$d\rho\co\mathfrak{g}\to\mathfrak{gl}(V)$ is faithful (Here
$\mathfrak{g}$ is the Lie algebra of $G$, $\mathfrak{gl}(V)$ the Lie
algebra of endomorphisms of $V$.
Let $\mathfrak{h}$ be a Cartan subalgebra of $\mathfrak{g}$ and let
$\Phi$ be the (reduced) root system relative to $\mathfrak{h}$.
Let $\Delta$ be a
choice of simple roots.  Then there is a \emph{Chevalley basis} 
for $\mathfrak{g}$ (see \cite{steinberg:notes} or \cite{chevalley:tohoku}) of
the form $\mathcal{B}=\{X_\alpha\ |\ \alpha\in\Phi\}\cup\{H_\alpha\ |\
\alpha\in\Delta\}$, so that the $H_\alpha$ generate $\mathfrak{h}$ (and
thus commute) and the structure constants are all integral.
In other words, the $\Z$--span of $\mathcal{B}$ is actually a Lie algebra over
$\Z$.
In particular, if $\{\beta - r_{\beta,\alpha} \alpha, \ldots \beta,
\ldots \beta + q_{\beta,\alpha} \alpha\}$ are all the roots on the
line $\{\beta + t\alpha\ |\ t\in\Z\}$, then:
\begin{enumerate}
\item $[X_\alpha, X_{-\alpha}] = H_\alpha$.
\item $[H_\alpha, X_\beta] = A(\alpha,\beta) X_\beta$, where
  $A(\alpha,\beta)$ is an integer determined by $\alpha$
  and $\beta$.
\item $[X_\alpha,X_\beta] = 0$ if $\alpha + \beta$ is not in $\Phi$;
  otherwise $[X_\alpha,X_\beta] = \pm (r_{\beta,\alpha}+1)X_{\alpha+\beta}$.
\end{enumerate}
Notice that the $X_\alpha$ are all ad-nilpotent, and that the 
$H_\alpha$ are ad-semisimple.

It can be shown (see, e.g. \cite[\S27]{humphreys:introduction})
that $V$ contains 
an ``admissible lattice'' $V_\Z$: This is a free $\Z$--module in $V$
which is invariant under  
$(d\rho(X_\alpha))^m / m!$ for any $\alpha$ and any $m$.  Let
$\{v_1,\ldots,v_n\}$ be a basis for $V_\Z$.  In terms of this basis,
$\rho(g)$ can be written as an $n$ by $n$ matrix with complex
entries for any $g$.  Let $x_{ij}\co G\to \C$ be the function which
simply reads off the $ij$th entry of this matrix.  The $x_{ij}$
generate an affine complex algebra $\C(G)$.  Let $\Z(G)$ denote the
$\Z$--algebra with the same generators.

We can endow $\Z(G)$ with a Hopf algebra structure by defining
a comultiplication $\mu^*$ by
$\mu^*(x_{ij}) = \sum_k x_{ik}\otimes x_{kj}$, a counit $\epsilon$
by $\epsilon(x_{ij}) = \delta_{ij}$, and an antipode $s$ by
$s(f)(g) = f(g^{-1})$.
(That $s$ maps  $\Z(G)$ into itself uses the fact
that $\det(\rho(g))=\pm 1$ for every $g\in G$ -- this 
 follows from semisimplicity.)
Since $\Z(G)$ is a Hopf algebra over $\Z$, it defines a functor from
rings to groups as follows:
For any ring $R$, let 
$G(R) = \mathrm{Hom}_\Z(\Z(G), R)$.  (Note that the elements of
$G(R)$ are $\Z$--algebra homomorphisms, so they send $1$ to $1$.)
We define a group operation
$\bullet$ on $G(R)$ by $\rho\bullet\sigma = (\rho\otimes\sigma)\circ
\mu^*$.
(In other words, $\rho\bullet\sigma(x_{ij}) = \sum_k
\rho(x_{ik})\sigma(x_{kj})$.) 

Some observations:
\begin{enumerate}
\item $G(\C)$ can be naturally identified with the image of $\rho\co
  G\to GL(V)$ -- if $\rho$ is assumed to be faithful, then clearly
  $G(\C)\cong G$ doesn't depend on $\rho$.  On the other hand, for $R$
  arbitrary, $G(R)$ does  depend  on $\rho$.
\item  The assignment $R\mapsto G(R)$ is a covariant functor from
  commutative rings to groups.  This functor is often called a
  ``Chevalley-Demazure group scheme.''
\item\label{zoft}  Any Hopf algebra gives such a functor.  An
  important example is
  $\Z[t]$, with a Hopf ($\Z$--)algebra structure so that:
\begin{enumerate}
\item The comultiplication satisfies $\mu^*(t)=t\otimes 1 + 1\otimes
  t$.
\item The counit satisfies $\epsilon(t)=0$.
\item The antipode satisfies $S(t)=-t$
\end{enumerate}
In this case the functor $R\mapsto \mathrm{Hom}(Z[t],R)$ is the
  ``forgetful'' functor, which takes a ring to its underlying Abelian
  group.
\item\label{contrafunctor} Since $G(R)$ is $Hom_\Z(\Z(G),R)$, any morphism of Hopf algebras
 $\Z(G)\to A$ where $A$ is
some other Hopf algebra, will give rise to a homomorphism of groups
 $Hom_\Z(A,R)\to G(R)$. 
\end{enumerate}
\begin{definition}
Let $\alpha$ be a root.  Then since $d\rho(X_\alpha)$ is
nilpotent, the formal sum $\exp ( t X_\alpha ) = \sum_{m=1}^\infty
t^m \frac{d\rho(X_\alpha)^m}{m!} $ is a matrix with entries which are
polynomials in $t$.  
Because $\frac{d\rho(X_\alpha)^m}{m!}$ preserves
$V_\Z$, each $x_{ij}(\frac{d\rho(X_\alpha)^m}{m!})$ is an integer.
Thus we get a map
\[\mathrm{ev}_\alpha\co \Z(G) \to \Z[t],\]
which sends $x_{ij}$ to the (integral) polynomial in $t$ which appears
in the $ij$'th place of the matrix $\exp( t X_\alpha)$.  The map
$\mathrm{ev}_\alpha$ is a morphism of Hopf algebras, where $\Z[t]$ is given a
Hopf algebra structure as in observation \eqref{zoft} above.
This morphism gives rise via observation \eqref{contrafunctor} to a
homomorphism of the additive group underlying $R$ into the group
$G(R)$
\[x_\alpha\co R\to G(R).\]  
The image of this map is the \emph{root subgroup}
of $G(R)$ corresponding to $\alpha$.
\end{definition}

\begin{definition}
The subgroup of $G(R)$ generated by the root subgroups
is denoted $E(R)$.  
\end{definition}

In this paper, we focus mainly on the special case that $R$ is the
ring of integers of a number field $k$ and $\Phi$ has rank at least
two.  In this case, $G(R)=E(R)$ by 
a result of Matsumoto \cite{matsumoto:csp}.  By a result 
of Carter and Keller \cite{carterkeller:slno} in case $\Phi=A_n$,
and  Tavgen$'$ \cite{tavgen:bg} in general,
$G(R)$ is boundedly generated by the root subgroups.

It should also be noted that in this special case, $G(R)$ is actually
an irreducible lattice in a semisimple Lie group.  Indeed, if $s$ and
$t$ are the number of real 
and complex places of $k$, respectively, then $G(R)$ is a lattice in the
Lie group 
\[ G(\R)^s\times G(\C)^t.\]

\section{Equivalence of actions}\label{s:actions}
We wish to study the variety of actions of a group $G$ on
Gromov hyperbolic spaces, up to some kind of coarse equivalence.  By a
{\em (Gromov) hyperbolic $G$--space}, we will always mean a Gromov
hyperbolic geodesic metric space, equipped with an isometric
$G$--action.

\begin{definition}
Let $X$ and $Y$ be Gromov hyperbolic $G$--spaces.  We say that $X$ and
$Y$ are {\em equivalent} if they lie in the same
equivalence class, under the equivalence relation generated by
coarsely equivariant quasi-isometric embeddings.
\end{definition}
The following proposition should serve to clarify this equivalence
relation. 
\begin{proposition}If $X_1$ and $X_2$ are equivalent Gromov hyperbolic
  $G$--spaces, then there is
  a third Gromov hyperbolic $G$--space $V$ which coarsely equivariantly
  quasi-isometrically embeds in both $X_1$ and $X_2$.
\end{proposition}
\begin{proof}
The key claim is the following:
\begin{claim}\label{c:diamond}
If $V$ and $W$ are hyperbolic $G$--spaces which coarsely equivariantly
  quasi-isometrically embed in a third hyperbolic $G$--space $X$, then
  there is a fourth hyperbolic $G$--space $A$ which coarsely
  equivariantly quasi-isometrically embeds into $V$ and $W$.
\end{claim}
Before proving the claim, let us see how it implies the proposition.
If $X_1$ and $X_2$ are equivalent hyperbolic $G$--spaces, they 
must be joined by a sequence of coarsely
equivariant quasi-isometric embeddings:
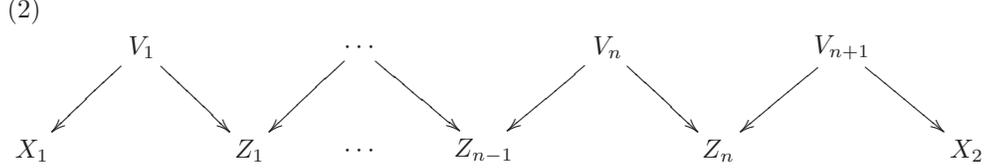
\cdlabel{long}{ & V_1\ar[dl]\ar[dr] & & \cdots\ar[dl]\ar[dr] & &V_n\ar[dl]\ar[dr] & &
  V_{n+1}\ar[dl]\ar[dr] & \\
X_1 & & Z_1 &\cdots & Z_{n-1} & & Z_n & & X_2}
By applying the claim with $V=V_n$, $W=V_{n+1}$, and $X=Z_n$, we
obtain a hyperbolic $G$--space $A$ which coarsely equivariantly
quasi-isometrically embeds into both $V_n$ and $V_{n+1}$, and hence
into both $Z_{n-1}$ and $X_2$.  We can thus shorten the sequence
\eqref{long}, unless it is a shortest possible such sequence,
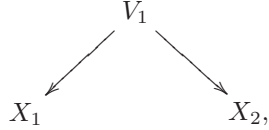
\cd{ & V_1\ar[dl]\ar[dr] & \\
 X_1 & & X_2,}
in which case the proposition is verified.

\begin{proof} (of Claim \ref{c:diamond})
We first construct a $G$--space $A_1$ which coarsely equivariantly
quasi-isometrically embeds in both $V$ and $W$ and is Gromov
hyperbolic but not geodesic.  We then show that $A_1$ is
quasi-isometric to a geodesic $G$--space $A$.

Choose $\delta>0$ so that all of $V$, $W$, and $X$ are
$\delta$--hyperbolic spaces.
By hypothesis, we may choose $K\geq 1$ and $C\geq 0$ so that there are maps
$\phi\co V\to X$ and $\psi\co W\to Y$ which are $C$--coarsely equivariant
$(K,C)$--quasi-isometric embeddings.  We choose constants $J_0<J_1<J_2$:
Let
$J_0 = 2B(K,C,\delta)+2\delta$,
where $B(K,C,\delta)$ is the constant of quasi-geodesic stability from
Proposition \ref{p:quasistable}, let $J_1=J_0+2C$, and let $J_2=4 J_1$.

Let $A_0$ be the subset of $V\times
W$ given by $A_0=\{(v,w)\mid d(\phi(v),\psi(w))\leq J_0\}$, and let $A_1$
be the smallest $G$--equivariant subset of $V\times W$ containing
$A_0$.  We endow $A_1$ with the (pseudo)metric
\[
 d((v_1,w_1),(v_2,w_2))=d(\phi(v_1),\phi(v_2))+d(\psi(v_2),\psi(w_2)).
\]
\begin{subclaim}
$A_1$ is Gromov hyperbolic.
\end{subclaim}
\begin{proof}
Since $A_1$ is not a geodesic space, we must work with the four-point
definition \ref{d:fourpoint}.  We will show that
$A_1$ is $(2\delta+4J_1)$--hyperbolic in the sense of Definition
\ref{d:fourpoint}.  
Let $\{p_i=(v_i,w_i)\mid i=1,\ldots,4\}$ be four points in $A_1$.  
For each $i$, we write $\uu{v}_i$ for $\phi(v_i)$ and $\uu{w}_i$ for
$\psi(w_i)$. By
reordering the points if necessary, we can assume that
\[d(\uu{v}_1,\uu{v}_3)+d(\uu{v}_2,\uu{v}_4)\geq
d(\uu{v}_1,\uu{v}_2)+d(\uu{v}_3,\uu{v}_4),\]
as in Figure \ref{f:palimpsest}.
\begin{figure}[htbp]
\begin{center}
\begin{picture}(0,0)%
\epsfig{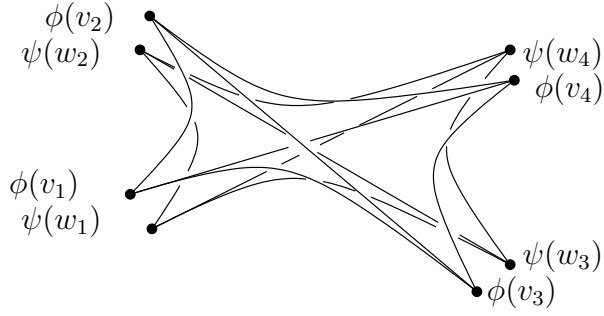}%
\end{picture}%
\setlength{\unitlength}{3947sp}%
\begingroup\makeatletter\ifx\SetFigFont\undefined%
\gdef\SetFigFont#1#2#3#4#5{%
  \reset@font\fontsize{#1}{#2pt}%
  \fontfamily{#3}\fontseries{#4}\fontshape{#5}%
  \selectfont}%
\fi\endgroup%
\begin{picture}(3300,1939)(2551,-2744)
\put(5551,-2686){\makebox(0,0)[lb]{\smash{{\SetFigFont{12}{14.4}{\familydefault}{\mddefault}{\updefault}{\color[rgb]{0,0,0}$\phi(v_3)$}%
}}}}
\put(2776,-961){\makebox(0,0)[lb]{\smash{{\SetFigFont{12}{14.4}{\familydefault}{\mddefault}{\updefault}{\color[rgb]{0,0,0}$\phi(v_2)$}%
}}}}
\put(2626,-1186){\makebox(0,0)[lb]{\smash{{\SetFigFont{12}{14.4}{\familydefault}{\mddefault}{\updefault}{\color[rgb]{0,0,0}$\psi(w_2)$}%
}}}}
\put(2551,-2011){\makebox(0,0)[lb]{\smash{{\SetFigFont{12}{14.4}{\familydefault}{\mddefault}{\updefault}{\color[rgb]{0,0,0}$\phi(v_1)$}%
}}}}
\put(2626,-2236){\makebox(0,0)[lb]{\smash{{\SetFigFont{12}{14.4}{\familydefault}{\mddefault}{\updefault}{\color[rgb]{0,0,0}$\psi(w_1)$}%
}}}}
\put(5851,-1411){\makebox(0,0)[lb]{\smash{{\SetFigFont{12}{14.4}{\familydefault}{\mddefault}{\updefault}{\color[rgb]{0,0,0}$\phi(v_4)$}%
}}}}
\put(5776,-1186){\makebox(0,0)[lb]{\smash{{\SetFigFont{12}{14.4}{\familydefault}{\mddefault}{\updefault}{\color[rgb]{0,0,0}$\psi(w_4)$}%
}}}}
\put(5776,-2461){\makebox(0,0)[lb]{\smash{{\SetFigFont{12}{14.4}{\familydefault}{\mddefault}{\updefault}{\color[rgb]{0,0,0}$\psi(w_3)$}%
}}}}
\end{picture}%
\caption{Points in $X$.}
\label{f:palimpsest}
\end{center}
\end{figure}
There are then two cases.

In case
\[ d(\uu{w}_1,\uu{w}_3)+d(\uu{w}_2,\uu{w}_4)\geq
d(\uu{w}_1,\uu{w}_2)+d(\uu{w}_3,\uu{w}_4),\]
then 
$\max\{d(p_1,p_3)+d(p_2,p_4),d(p_1,p_2)+d(p_3,p_4)\}=d(p_1,p_2)+d(p_2,p_4)$.
Applying $\delta$--hyperbolicity in $X$, we obtain
\begin{eqnarray*}
d(p_1,p_4)+d(p_3,p_2) & = &
d(\uu{v}_1,\uu{v}_4)+d(\uu{v}_2,\uu{v}_3)+d(\uu{w}_1,\uu{w}_4)+d(\uu{w}_2,\uu{w}_3)
\\
 & \leq & d(\uu{v}_1,\uu{v}_3)+d(\uu{v}_2,\uu{v}_4)+ 2\delta
+d(\uu{w}_1,\uu{w}_3)+d(\uu{w}_2,\uu{w}_4)+ 2\delta \\
& = & d(p_1,p_3)+d(p_2,p_4)+ 4\delta \\ 
& = & \max\{d(p_1,p_3)+d(p_2,p_4),d(p_1,p_2)+d(p_3,p_4)\} + 2
(2\delta) \\ 
& \leq & \max\{d(p_1,p_3)+d(p_2,p_4),d(p_1,p_2)+d(p_3,p_4)\} \\
& & + 2(2\delta+4 J_1)
\end{eqnarray*}
as required.

If on the other hand 
\[ d(\uu{w}_1,\uu{w}_3)+d(\uu{w}_2,\uu{w}_4) <
d(\uu{w}_1,\uu{w}_2)+d(\uu{w}_3,\uu{w}_4),\]
then it is still true (since $d(\uu{w_i},\uu{v_i})\leq J_1)$) that 
\begin{eqnarray*}
d(\uu{w}_1,\uu{w}_3)+d(\uu{w}_2,\uu{w}_4) & \geq &
d(\uu{v}_1,\uu{v}_3)+d(\uu{v}_2,\uu{v}_4) - 4J_1 \\
& \geq & d(\uu{v}_1,\uu{v}_2)+d(\uu{v}_3,\uu{v}_4) - 4J_1 \\
& \geq & d(\uu{w}_1,\uu{w}_2)+d(\uu{w}_3,\uu{w}_4) - 8J_1.
\end{eqnarray*}
We therefore obtain
\begin{eqnarray*}
d(p_1,p_4)+d(p_3,p_2) & = &
d(\uu{v}_1,\uu{v}_4)+d(\uu{v}_2,\uu{v}_3)+d(\uu{w}_1,\uu{w}_4)+d(\uu{w}_2,\uu{w}_3)
\\
& \leq & d(\uu{v}_1,\uu{v}_3)+d(\uu{v}_2,\uu{v}_4)+ 2\delta
+d(\uu{w}_1,\uu{w}_2)+d(\uu{w}_3,\uu{w}_4)+ 2\delta \\
& \leq & d(p_1,p_3)+d(p_2,p_4) + 2(2\delta+4 J_1) \\
& \leq & \max\{d(p_1,p_3)+d(p_2,p_4),d(p_1,p_2)+d(p_3,p_4)\}\\& & + 2
(2\delta+4 J_1),
\end{eqnarray*}
which finishes the proof of the Subclaim.
\end{proof}
Though the space $A_1$ is hyperbolic, and $G$ quasi-acts on $A_1$
via the diagonal action on $V\times W$, the space $A_1$ is not  geodesic,
and $G$ does not act by isometries.  Both of these issues can be fixed
at once, by replacing $A_1$ by an appropriate graph.  
Specifically, we let $A$ be the graph with vertex set
$V(A)=A_1$, and with an edge between every pair of vertices $p$ and
$q$ so that there exists some $g$ with $d(gp,gq)\leq J_2:=4J_1$.  Clearly
this graph is a geodesic $G$--space.

Let $\iota\co A_1\to A$ be the map which takes a point to the
corresponding vertex.  We claim that $\iota$ is a quasi-isometry, and
so $A$ is a Gromov hyperbolic space.

Let $a=(v_a,w_a)$ and $b=(v_b,w_b)$ be points in $A_1$ so that
$d(\iota(a),\iota(b))=1$.  There is some $g$ with 
\[d(ga,gb)=d(\phi(gv_a),\phi(gv_b))+d(\psi(gv_a),\psi(gv_b))\leq J_2.\]
Since $\phi$ and $\psi$ are both $C$--coarsely equivariant,
$d(\phi(gv_a),\phi(gv_b))$ differs  by
at most $2C$ from $d(\phi(v_a),\phi(v_b))$; similarly, $d(\psi(gw_a),\psi(gw_b))$ differs from $d(\psi(w_a),\psi(w_b))$  by
at most $2C$.  It follows that $d(ga,gb)$ differs by at most $4C$ from
$d(a,b)$, and so $d(a,b)\leq J_2+4C<\frac{3}{2}J_2$.  As this is true for every pair
of points connected by an edge, we deduce
\begin{equation}\label{lower}
 d(\iota(p),\iota(q)) \geq \frac{2}{3J_2} d(p,q) 
\end{equation}
for any pair of points $p$, $q\in A_1$.

We now obtain the complementary bound to \eqref{lower}.  Suppose $p=(v_1,w_1)$
and $q=(v_2,w_2)$ are any two points in $A_1$.  We write $\uu{v_1}$
for $\phi(v_1)$ and so on as before.
We will show that
\begin{equation}\label{upper}
d(\iota(p),\iota(q))\leq  \frac{1}{J_1}d(p,q)+2
\end{equation}
by constructing a path joining $\iota(p)$ to $\iota(q)$ in $A$.  The
vertices of this path will be points in $A_0\subset V\times W$ so that the first
coordinate lies on a geodesic between $v_1$ and $v_2$, while the
second lies on a geodesic between $w_1$ and $w_2$.

If $d(p,q)\leq 2J_2 = 8J_1$, then \eqref{upper} is automatically
satisfied.  We may therefore assume that $d(p,q)\geq 2J_2$.  Because
$d(\uu{v}_1,\uu{w}_1)$ and $d(\uu{w}_1,\uu{w}_2)$ are at most $J_0$,
we have
$\min\{d(\uu{v}_1,\uu{v}_2),d(\uu{w}_1,\uu{w}_2)\}\geq J_2-J_0 >
3J_0$. It follows that
$\gp{\uu{v}_2}{\uu{w}_2}{\uu{v}_1}>\gp{\uu{w}_1}{\uu{w}_2}{\uu{v}_1}$, as in
Figure \ref{f:bar}.
\begin{figure}[htbp]
\begin{center}
\input{bar.pstex_t}
\caption{A pair of pairs of points in $X$, together with comparison tripods.
The dashed line on the left is the image of a geodesic in $V$; the one
on the right is the image of a geodesic in $W$.}
\label{f:bar}
\end{center}
\end{figure}
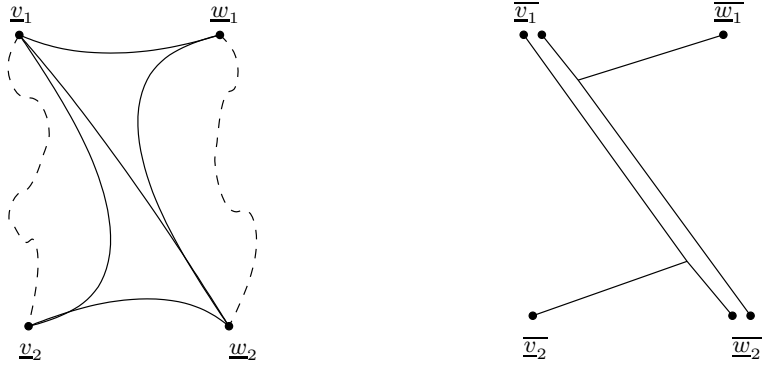
In fact,
$\gp{\uu{v}_2}{\uu{w}_2}{\uu{v}_1}-\gp{\uu{w}_1}{\uu{w}_2}{\uu{v}_1}>J_0\geq
J_1/2$ (since clearly $B(K,C,\delta)\geq C$).
We may thus choose real numbers 
\[ \gp{\uu{w}_1}{\uu{w}_2}{\uu{v}_1}=t_0<t_1<\cdots<t_k =
\gp{\uu{v}_2}{\uu{w}_2}{\uu{v}_1}\]
satisfying $\frac{J_1}{2}<|t_{i+1}-t_i|\leq J_1$ for each $i$.
Note that $k \leq \frac{1}{J_1}d(p,q)$ for such a choice.  

Let
$\gamma$ be a unit speed geodesic from $\uu{v}_1$ to $\uu{w}_2$.  For
each $i$ between $0$ and $k$, there are points $x_i'$ on a geodesic
between $\uu{v}_1$ and $\uu{w}_1$ and $y_i'$ on a geodesic between
$\uu{w}_1$ and $\uu{w}_2$ which are distance at most $\delta$ from
$\gamma(t_i)$.  Since $V$ and $W$ are geodesic spaces, there are
geodesics $[v_1,v_2]$ in $V$ and $[w_1,w_2]$ in $W$.  The maps $\phi$
and $\psi$ send these geodesics
to $(K,C)$--quasi-geodesics in $X$.  Applying quasi-geodesic stability
(Proposition \ref{p:quasistable}), there are points
$x_i\in [v_1,v_2]$, $y_i\in [w_1,w_2]$ so that $\phi(x_i)$ and
$\psi(y_i)$ lie within $B+\delta$ of $\gamma(t_i)$.  Since
$d(\phi(x_i),\psi(y_i))\leq 2B+2\delta=J_0$, the point $p_i=(x_i,y_i)$
lies in $A_0\subseteq A_1$.  Moreover, $d(p_i,p_{i+1})\leq
2(2(B+\delta)+J_1)<J_2$, and so $d(\iota(p_i),\iota(p_{i+1}))\leq 1$.

Finally, one notes that 
\[
d(\uu{v}_1,\phi(x_0))+d(\uu{w}_1,\psi(y_0))  \leq  J_0+2\delta+2 B=2J_0<J_2
\]
and likewise for $d(\uu{v}_1,\phi(x_0))+d(\uu{w}_1,\psi(y_0))$, from
which it follows that $d(\iota(p),\iota(p_0))$ and
$d(\iota(q),\iota(p_k))$ are at most one.  Thus
\[ d(\iota(p),\iota(q))\leq k+2 \leq \frac{1}{J_1} d(p,q)+2.\]
It is obvious that $\iota$ is $1$--almost onto, and so $\iota$ is a
quasi-isometry.

It follows that $A$ is a Gromov hyperbolic geodesic $G$--space
which coarsely equivariantly, quasi-isometrically embeds into $V$ and
$W$.
  
\end{proof}
\end{proof}

\begin{remark}
We record some observations about this equivalence:
\begin{enumerate}
\item Let $\Gamma$ be a  group acting isometrically on $\H^2$.  This
  action extends in an obvious way to either $\H^3$ or $\C\H^2$.
  Although there is no quasi-isometric embedding either from $\H^3$ to
  $\C\H^2$ or vice versa, these actions are equivalent under the
  equivalence relation.
\item A Gromov hyperbolic $G$--space has a bounded orbit if and only if
  it is equivalent to the trivial $G$--space consisting of a single
  point.
\item If $X$ and $Y$ are equivalent Gromov hyperbolic $G$--spaces,
  $x\in X$ and $y\in Y$, then the
  limit sets $\Lambda(X)=\{e\in \partial X\cap \overline{Gx}\}$ and
  $\Lambda(Y)=\{e\in \partial Y\cap \overline{Gy}\}$ are equivariantly
  homeomorphic.
\item The equivalence is perfectly well-defined in the more general
  setting of {\em quasi-actions} on geodesic Gromov hyperbolic
  spaces.  Call a geodesic Gromov hyperbolic space with a
  $G$--quasi-action a {\em hyperbolic quasi-$G$--space}.  Every
  hyperbolic quasi-$G$--space is equivalent to some hyperbolic
  $G$--space. 
\end{enumerate}
\end{remark}

Given a group acting on a hyperbolic $G$--space $X$ and some basepoint
$x_0\in X$, one obtains a metric on $G$ given by the formula
\[d(g,h) = d(x_0,g^{-1}hx_0).\]
This metric is obviously determined by its values on $\{1\}\times G$.
The compact-open topology on real-valued functions on $G$ thus induces
a topology on the set of pointed hyperbolic $G$--spaces.
The quotient topology on the space of equivalence classes
of Gromov hyperbolic $G$--spaces is not Hausdorff.   For example, if
$G$ is a surface 
group, then all of Teichm\"uller space is identified to a single
point, whose closure contains many inequivalent actions of $G$ on
$\R$--trees.

\subsection{Combinatorial horoballs}\label{ss:combhoro}
Combinatorial horoballs of the simplest possible type were defined in
\cite{rhds}, and used as building blocks for complexes naturally
associated to relatively hyperbolic groups.  The point there as here
is that these spaces can be used to ``hide'' an action on a
non-hyperbolic space in an action on a hyperbolic space.  There is
some flexibility as to how this can 
happen which is deliberately ignored in \cite{rhds}; here we give a
more general construction.  

\begin{definition}\label{d:admissible}
Let $X$ be a graph, acted on by $G$, and suppose that 
\[\{B_i\co  X^{(0)}\to 2^{X^{(0)}}\}_{i\in \N}\]
is a collection of functions.  We will say that $B_*$ is
\emph{admissible} if it
satisfies the following
four axioms:
\begin{enumerate}
\item Connectedness:  If $v\in X^{(0)}$, then
  $B_1(v)=\{w\in X^{(0)}\mid d_X(v,w)\leq 1\}$.
\item\label{ax:exp} Exponential growth: Let $v$, $w\in X^{(0)}$ and $n\in \N$.  If
  $v\in B_n(w)$, then $B_n(v)\subset  B_{n+1}(w)$. 
\item Symmetry: Let $v$, $w\in X^{(0)}$ and $n\in \N$.
  If $v\in B_n(w)$, then $w\in B_n(v)$.
\item $G$--equivariance: If $w\in X^{(0)}$, $n\in \N$, and $g\in G$,
  then $g B_n(w)=B_n(gw)$.
\end{enumerate}
\end{definition}

\begin{definition}\label{d:choroball}
Let $B_*$ be a sequence of functions $\{B_i\co  X^{(0)}\to
2^{X^{(0)}}\}_{i\in \N}$.  The \emph{combinatorial horoball based on
  $X$ and $B_*$}, or $\mc{H}(X,B_*)$, is the graph defined as follows:
\begin{enumerate}
\item $\mc{H}(X,B_*)^{(0)} = X^{(0)}\times \N$.
\item If $n\in \N$ and $v\in X^{(0)}$, then $(v,n)$ is connected to
  $(v,n+1)$ by an edge (called a \emph{vertical edge}).
\item If $n\in \N$ and $v\in B_n(w)$, then $(v,n)$ is connected to
  $(w,n)$ by an edge (called a \emph{horizontal edge}).
\end{enumerate}
\end{definition}
We leave the following Lemma as an exercise.
\begin{lemma}\label{l:horohyp}
Let $G$ act on the graph $X$.  If $B_*$ is admissible, then
$\mc{H}(X,B_*)$ is Gromov hyperbolic, and the action of $G$ on $X$
induces an action of $G$ on $\mc{H}(X,B_*)$.
\end{lemma}

Because of the flexibility of this construction, a group
typically admits many inequivalent actions on horoballs.

\begin{definition}\label{d:ihoroball}
A Gromov hyperbolic $G$--space {\em has an invariant horoball} if it is
equivalent to a $G$--space of the form $\mc{H}(X,B_*)$, for some graph
$X$ with a $G$--action, and some family $B_*$ which is admissible.
\end{definition} 
In the next section we give some other characterizations of $G$--spaces
with invariant horoballs (Theorem \ref{t:horeq}).

\section{Elementary actions and pseudocharacters}\label{s:elementary}
In this section we show how an elementary action by $G$ on a hyperbolic space
gives rise to a pseudocharacter on $G$ which ``picks out'' the
elements which act hyperbolically.  We then see that if no element
acts hyperbolically, then the action has an invariant horoball.

\subsection{The pseudocharacter coming from an elementary action}
Recall the following definitions.
\begin{definition}
An action of $G$ on a hyperbolic space $X$ is \emph{elementary} if it
is either equivalent to the trivial action on a point, or if the
induced action on $\partial X$ has a global fixed point.\footnote{This definition is slightly more restrictive than the
  usual one, which allows for a pair of points in $\partial X$ to be
  preserved.} 
\end{definition}
\begin{definition}
A \emph{quasicharacter} (or quasi-morphism) on a group $G$ is a
real valued function $q$ on $G$ satisfying 
\begin{equation}\label{defect}
|q(gh)-q(g)-q(h)|<C\mbox{, for all } g,h\in G.
\end{equation}
The \emph{defect} of a quasicharacter is the
smallest $C$ so that \eqref{defect} is satisfied.  A
\emph{pseudocharacter} (or homogeneous quasi-morphism) is a
quasicharacter $p$ which satisfies the additional condition
\begin{equation*}
  p(g^n)=n p(g)\mbox{, for all }g\in G, n\in \Z
\end{equation*}
\end{definition}
In this subsection we show how an elementary action on a hyperbolic
space gives rise to a pseudocharacter which is nonzero precisely on
the elements which act hyperbolically.
We begin by studying ``quasi-horofunctions'' on the space
$X$, corresponding to a fixed point at infinity.  A
quasi-horofunction restricted to an arbitrary orbit will give 
a quasicharacter, which can then be homogenized to give the desired
pseudocharacter.
  
\begin{definition}
(cf. \cite[7.5.D]{gromov:wordhyperbolic})
Let $\mathbf{x}=\{x_i\}$ be a sequence tending to infinity in the geodesic
hyperbolic space $X$.  The \emph{quasi-horofunction coming from}
$\mathbf{x}$ is the function $\eta_{\mathbf{x}}\co X\to \R$ given by
\[ \eta_{\mathbf{x}}(a) = \limsup_{n\to\infty}(d(a,x_n)-d(x_0,x_n)). \]
\end{definition}

We use the following observation repeatedly:
\begin{observation}\label{o:fourdelta}
Let $A$, $B$, $C$ and $D$ be four points in the $\delta$--hyperbolic
space $X$. If $\gp{C}{D}{A}$ and $\gp{C}{D}{B}$ are both larger than $d(A,B)$,
then 
\[ | (d(B,C) - d(A,C)) - (d(B,D) - d(A,D)) | \leq 4\delta. \]
\end{observation}
The observation \ref{o:fourdelta} implies in particular:
\begin{lemma}\label{l:approxbyterm}
If $a\in X$, and $\mathbf{x}=\{x_i\}$ tends to infinity in $X$, then
for all $n$ sufficiently large, 
\[| \eta_{\mathbf{x}}(a)-(d(a,x_n)-d(x_0,x_n))| \leq 4\delta.\]
\end{lemma}

We now can describe the dependence of $\eta_{\mathbf{x}}$ on the sequence
$\mathbf{x}$.
\begin{lemma}\label{l:key}
Let $\mathbf{x}=\{x_i\}$ and $\mathbf{y}=\{y_i\}$ be two sequences of
points in the geodesic $\delta$--hyperbolic space $X$
which tend to the same point in $\partial X$.  For any point $a\in X$,
we have 
\[ | \eta_{\mathbf{x}}(a)-\eta_{\mathbf{y}}(a)-\eta_{\mathbf{x}}(y_0)|
\leq 16\delta. \]
\end{lemma}
\begin{proof}
Since $\mathbf{x}$ and $\mathbf{y}$ tend to the same point at
infinity, we may choose $N$ so that 
$(z,z')_\alpha> 2\thinspace \mathrm{diam}\{a,x_0,y_0\}$ for every $z$, $z'$ in
$\{x_i \mid i\geq N\}\cup \{y_i \mid i\geq N\}$.  Using Lemma
\ref{l:approxbyterm} three times, the quantity 
\[|\eta_{\mathbf{x}}(a)-\eta_{\mathbf{y}}(a)-\eta_{\mathbf{x}}(y_0)|\]
differs by at most $12\delta$ from
\begin{equation}\label{easy}
| (d(a,x_N) - d(a,y_N)) - (d(y_0,x_N)-d(y_0,y_N))|.
\end{equation}
By Observation \ref{o:fourdelta}, the quantity \eqref{easy} is at most
$4\delta$.  The Lemma follows.
\end{proof}

Using Lemma \ref{l:key}, we deduce that an isometry of $X$ changes
$\eta_{\mathbf{x}}(a)$ by approximately the same amount, independent of
the $a\in X$ chosen:
\begin{proposition}\label{p:isom}
Let $X$ be a geodesic $\delta$--hyperbolic space, and suppose that
$\mathbf{x}=\{x_i\}$ tends to $e\in \partial X$.  Let $a$ be any point
in $X$.
If $g\in \mathrm{Isom}(X)$ fixes $e$, then
$\eta_{\mathbf{x}}(ga)$ differs from
$\eta_{\mathbf{x}}(a)+\eta_{\mathbf{x}}(gx_0)$ by at most $16\delta$.
\end{proposition}\label{p:trans}
\begin{proof}
First note that if $g\mathbf{x}$ is the sequence $\{gx_i\}$, then
\[\eta_{\mathbf{x}}(a) = \eta_{g\mathbf{x}}(ga).\]
But by Lemma \ref{l:key}, 
\[ | \eta_{\mathbf{x}}(ga)-\eta_{g\mathbf{x}}(ga)-\eta_{\mathbf{x}}(g x_0)|
\leq 16\delta. \]
\end{proof}

\begin{corollary}\label{c:quasimorphism}
Suppose $X$ is a $\delta$--hyperbolic space, and that $G$ acts on $X$
fixing $e\in \partial X$.  Let $\mathbf{x}=\{x_i\}$ be any sequence tending to
$e$.  The function $q_{\mathbf{x}}\co G\to \R$ defined by
$q_{\mathbf{x}}(g)=\eta_{\mathbf{x}}(g x_0)$ is a quasicharacter of
defect at most $16\delta$.
\end{corollary}
\begin{proof}
Let $g$, $h\in G$.  Using Proposition \ref{p:isom},
\begin{eqnarray*}
|q_{\mathbf{x}}(gh)-q_{\mathbf{x}}(g)-q_{\mathbf{x}}(h)| & = &
 |\eta_{\mathbf{x}}(ghx_0)-\eta_{\mathbf{x}}(gx_0)-\eta_{\mathbf{x}}(hx_0)| \\
& \leq & |
 \eta_{\mathbf{x}}(hx_0)+\eta_{\mathbf{x}}(gx_0)-\eta_{\mathbf{x}}(gx_0)-\eta_{\mathbf{x}}(hx_0)|+16\delta = 16\delta.
\end{eqnarray*}
\end{proof}

\begin{proposition}\label{p:pseudo}
Let $X$, $e$, $\mathbf{x}$, and $q_{\mathbf{x}}$ be as in the
statement of Corollary \ref{c:quasimorphism}, and let the
pseudocharacter $p_{\mathbf{x}}\co G\to \R$
be given by
\[p_{\mathbf{x}}(g) = \lim_{n\to\infty}\frac{q_{\mathbf{x}}(g^n)}{n}.\]
Then $p_{\mathbf{x}}(g)\neq 0$ if and only if $g$ acts hyperbolically on $X$.
\end{proposition}

\begin{proof}
First, we suppose that $p_{\mathbf{x}}(g)\neq 0$.  Without loss of generality we
assume that $p_{\mathbf{x}}(g)>0$.  Since
\[p_{\mathbf{x}}(g)=\lim_{n\to\infty}\frac{\eta_{\mathbf{x}}(g^n x_0)}{n}>0,\]
there exists some $N$ so that
$\eta_{\mathbf{x}}(g^n x_0)>\frac{1}{2}p_{\mathbf{x}}(g)n$ for all
$n\geq N$.

Let $a$ and $b$
be integers.  Choosing some sufficiently large $M$ and applying the triangle
inequality and Lemma \ref{l:approxbyterm}, we obtain a lower bound for
$d(g^a x_0,g^b x_0)$:
\begin{eqnarray*}
d(g^a x_0,g^b x_0) & = & d(g^N x_0,g^{N+|b-a|}x_0)\\
 & \geq & d(g^{N+|b-a|} x_0, x_M)-d(g^N x_0,x_M)\\
 & \geq & \eta_{\mathbf{x}}(g^{N+|b-a|})-\eta_{\mathbf{x}}(g^N) -
 8\delta\\
 & \geq & \frac{1}{2}p_{\mathbf{x}}(g)(N+|b-a|)-\eta_{\mathbf{x}}(g^N)-8\delta\\
 & = & \frac{1}{2}p_{\mathbf{x}}(g) |b-a| - (\eta_{\mathbf{x}}g^N+8\delta-\frac{1}{2}p_{\mathbf{x}}(g)N).
\end{eqnarray*}
On the other hand, $d(g^a x_0,g^b x_0)\leq |b-a| d(x_0,gx_0)$, so the map 
$n\mapsto g^n(x_0)$ is a quasi-isometric embedding, and $g$ acts
hyperbolically.

Conversely, suppose that $g$ acts hyperbolically.  
It
follows that there is some $\epsilon>0$ so that
\[d(g^n x_0,x_0)>\epsilon n\] for all $n$.
By
replacing $g$ with $g^{-1}$, we may suppose that $\{g^i x_0\}$ tends
to a point in $\partial X\smallsetminus\{e\}$ as $i\to\infty$.
Thus there is some $R$ so that 
\[\gp{g^n x_0}{x_i}{{x_0}}<R\] for all positive $n$ and $i$.
Lemma \ref{l:approxbyterm} implies that for sufficiently
large $m$,
\begin{eqnarray*}
\eta_{\mathbf{x}}(g^n x_0) & \geq & d(g^n x_0,x_m)-d(x_0,x_m)-4\delta\\
& = & d(g^n x_0,x_0)-2 \gp{g^nx_0}{x_n}{{x_0}} - 4\delta \\
& \geq & \epsilon n - (2 R +4\delta).
\end{eqnarray*}
Since $q_{\mathbf{x}}(g^n)\geq \epsilon n - (2R+4\delta)$ for all $n>0$, we must
have $p(g)\geq \epsilon>0$.
\end{proof}

\begin{remark}
Proposition \ref{p:pseudo} was proved in \cite[Proposition 3.9]{manning:qfa} for the case of
quasi-trees.  The proof here is somewhat more efficient even in this
case.
\end{remark}

\subsection{Characterization of $G$--spaces with invariant horoballs}
\begin{theorem}\label{t:horeq}
Let $X$ be a Gromov hyperbolic $G$--space.  The following are
equivalent:
\begin{enumerate}
\item\label{ih} $X$ has an invariant horoball.
\item\label{enh} $X$ is elementary, and no element acts
  hyperbolically.
\item\label{onepoint} $X$ is equivalent to a hyperbolic $G$--space $Y$, so that
  $\#(\partial Y)\leq 1$.
\end{enumerate}
\end{theorem}
\begin{proof}
That \eqref{ih} implies \eqref{onepoint} is trivial.

We next assume \eqref{onepoint} and show \eqref{enh}.  If \eqref{enh}
holds for $Y$, it holds for $X$, so we may suppose that $X=Y$.  If
$\#(\partial X) = 1$ or $G$ has a bounded orbit in $X$, then clearly
$X$ is elementary.  The only case remaining is that $\partial X$ is
empty, but $Gx$ is unbounded for some $x\in X$.  We show that this
case does not occur.  Chose a sequence
$\{g_i\}$ in $G$ so that $\lim_{i\to\infty}d(g_i x,x) = \infty$.
Since $\partial X$ is empty,
$\liminf_{i,j\to\infty}(g_ix,g_jx)_x<\infty$.  It follows that there
are elements $g_m$, $g_n$ so that $d(g_mx,x)$ and $d(g_nx,x)$ are much
larger than $(g_mx,h_nx)_x$.  It can then be shown (see, for
example, \cite[Chapitre 9, Lemme 2.3]{cdp}) that $g_mg_n$ is hyperbolic.
It follows that $\partial X$ contains at least two points (the fixed
points of $g_mg_n$), contrary to assumption.

It remains to show that \eqref{enh} implies \eqref{ih}.
Let $X$ be a hyperbolic $G$--space so that the action of $G$ is
elementary.  If $X$ is equivalent to a point (i.e. if $Gx$ is bounded
for $x\in X$), then $X$ is also
equivalent to a ray, which is the combinatorial horoball based on a
point.  We therefore may assume that $Gx$ is unbounded for any $x\in
X$.  We will construct a combinatorial horoball which coarsely
equivariantly quasi-isometrically embeds in $X$.
Let $\delta>0$ be some number so
that $X$ is $\delta$--hyperbolic.
Let $e$, $\mathbf{x}$,
$q_{\mathbf{x}}$ and $p_{\mathbf{x}}$ be as in the statements of
Corollary \ref{c:quasimorphism} and Proposition \ref{p:pseudo}.

To build the combinatorial horoball, we first must start with a graph $Y$
on which $G$ acts.  Choose a finite generating set $S$ for $G$, and
let $C_0 = \mathrm{diam}(S x_0)$. 
Let $V(Y)=G$, and connect $g$ to $h$ in $Y$ if $d(g x_0, h x_0)\leq C_0$.  It
is clear that $G$ acts on $Y$; in fact, $Y$ is a certain Cayley graph
for $G$.  We next define the functions $B_n\co V(Y)\to 2^{V(Y)}$.
Let $C_1 = 2C_0+20\delta$, 
and let
\begin{equation}\label{bstar} 
B_n(g) = \{ h\in G\mid d(hx_0,gx_0)\leq (2n+1)C_1\}.
\end{equation}
\begin{claim}
The sequence of functions $B_*$ in \eqref{bstar} is admissible in the sense
of Definition \ref{d:admissible}.
\end{claim}
\begin{proof}
The only axiom which is not obvious is \eqref{ax:exp}.  We must show
that if $a$ and $b$ are in $B_n(v)$, then $a\in B_{n+1}(b)$ (or
equivalently $b\in B_{n+1}(a)$).  
Put another way,
we must show that if $d(ax_0,vx_0)$ and $d(bx_0,vx_0)$
are bounded above by $(2n+1)C_1$ then $d(a x_0,b x_0)\leq (2n+3)C_1$.

Because no element of $G$ acts hyperbolically, the pseudocharacter
$p_{\mathbf{x}}$ is identically zero.  An easy argument shows that
$|q_{\mathbf{x}}(g)|\leq 16\delta$ for all $g\in G$.
Using Lemma \ref{l:approxbyterm}, we can choose some large $n$ so that 
\[| \eta_{\mathbf{x}}(z)-(d(z,x_n)-d(x_0,x_n))| \leq 4\delta\]
for $z\in \{ax_0,bx_0,vx_0\}$.  It follows that 
\begin{equation}\label{diff}
|d(z_1,x_m)-d(z_2,x_m)|\leq 16\delta+4\delta\leq C_1
\end{equation} for $z_1$,
$z_2\in\{ax_0,bx_0,vx_0\}$.  The assertion to be proved is symmetric
in $a$ and $b$, so we may assume
that $\gp{x_m}{b}{v}\leq \gp{x_m}{a}{v}$.  
We deduce:
\begin{eqnarray*}
d(a,b) & \leq & \gp{v}{x_m}{a}+[\gp{a}{x_m}{v}-\gp{b}{x_m}{v}]+\gp{v}{x_m}{b}+2\delta\\
& = & d(v,a) +(d(b,x_m)-d(v,x_m))+2\delta \\
& \leq & (2n+1)C_1 + C_1 + 2\delta \leq (2n+3)C_1.
\end{eqnarray*}
The first line follows from examining the comparison tripods
for the triangles $\Delta(a,v,x_m)$ and $\Delta(b,v,x_m)$;
the last follows from \eqref{diff}.
\end{proof}

Since $B_*$ is admissible, the combinatorial horoball
$H=\mathcal{H}_{B_*}(Y)$ is a hyperbolic $G$--space.  It remains
to construct a coarsely equivariant quasi-isometric embedding from $H$
to $X$.  It 
suffices to define this map on the vertices of $H$.
For each $g\in G = V(Y)$ and each $n\in\N$, choose some
$i(g,n)$ so that $\gp{x_k}{x_l}{{gx_0}}\geq 2nC_1$ for all $k$, $l\geq
i(g,n)$.  Choose also some unit speed geodesic $\gamma_{g,n}$ starting
at $g x_0$ and ending at $x_{i(g,n)}$.  Any vertex of $H$ is a pair
$(g,n)$ where $n\in \N$ and $g\in G$.  We define a map $\phi\co
V(H)\to X$ by
\begin{equation*}
\phi(g,n) = \gamma_{g,n}(nC_1).
\end{equation*}
A number of choices were made in the definition of $\phi$ (namely, the
sequence $\mathbf{x}$, the numbers $i(g,n)$, and the geodesics $\gamma_{g,n}$).
However, so long as $x_0$ is unchanged, different
choices lead to a function which differs by at most $\delta$
from $\phi$.  In particular, we could replace $\mathbf{x}$ by
$\mathbf{x}'=\{x_i'\}$, where $x_0'=x_0$ and $x_i'=hx_i$ for some
fixed $h$ and for all $i\geq 1$.
It follows that the distance between $\phi(hg,n)$ and $h\phi(g,n)$ is
at most $\delta$ for any $h$, $g\in G$ and $n\in \N$, and so the map
$\phi$ is coarsely equivariant.

It remains to show that $\phi$ is a quasi-isometric
embedding.

Note that if $v$ and $w$ are two vertices in $H$ connected by a
vertical path, then 
$C_1d_H(v,w)-\delta\leq d(\phi(v),\phi(w))\leq C_1d_H(v,w)+\delta$,
where $d_H$ is the distance in $H$. 

We therefore assume that $v=(a,n)$ and $w=(b,m)$, where $a\neq b$.
There is a unique $k$ so that $(2k-1)C_1<d(a,b)\leq (2k+1)C_1$.  If
$\max\{m,n\}\geq k$, then $d(v,w)=|m-n|+1$; otherwise
$d(v,w)=2k-(m+n)+\frac{1}{2}\pm\frac{1}{2}$.  Let 
$I \geq \max\{i(a,n),i(b,m)\}$, and observe that the points $\phi(v)$
and $\phi(w)$ lie within $\delta$ of geodesics joining $ax_0$ to $x_I$
and $bx_0$ to $x_I$, respectively (see Figure \ref{f:deep}).  
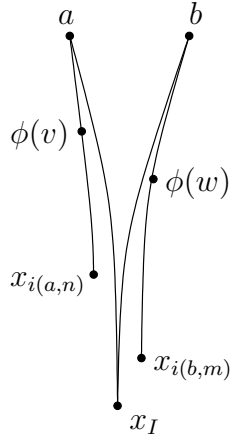
\begin{figure}[htbp]
\begin{center}
\input{deep.pstex_t}
\caption{A pair of points in the image of $\phi$.}
\label{f:deep}
\end{center}
\end{figure}
Note that $\gp{x_I}{ax_0}{{bx_0}}$
and $\gp{x_I}{bx_0}{{ax_0}}$ differ from $\frac{1}{2}d(ax_0,bx_0)$ by at
most $\frac{C_1}{2}$.
There are a couple of cases to consider.

First, assume that one or both of $n$ and $m$ is at least $k$.
Without loss of generality assume that $n\geq k$.  Since 
$d(ax_0,bx_0)\leq (2k+1)C_1$, we have $(x_I,bx_0)_a\leq kC_1+C_1$. 
Since $d(\phi(v),ax_0)\geq k C_1$,
it follows that $\phi(v)$ is at most $C_1+\delta$ from the
geodesic joining $bx_0$ to $x_I$. Accounting for the possible
difference between $\gp{x_I}{ax_0}{{bx_0}}$ and $\gp{x_I}{bx_0}{{ax_0}}$, we deduce that the distance between $\phi(v)$ and
$\phi(w)$ differs from $|n-m|C_1=(d(v,w)+\frac{1}{2}\pm\frac{1}{2})C_1$ by at most $2C_1+\delta$.  It
follows that in case one of $n$ or $m$ is at least $k$, we have
\[ C_1 d(v,w)-2C_1+\delta \leq d(\phi(v),\phi(w))\leq
C_1d(v,w)+3C_1+\delta.\]

In case both $n$ and $m$ are strictly less than $k$, we may argue as
follows.  Since $\gp{ax_0}{x_I}{{bx_0}}$ and $\gp{bx_0}{x_I}{{ax_0}}$ are both
at least $\frac{1}{2}(d(ax_0,bx_0)-C_1)\geq (k-1)C_1\geq \max\{m,n\}
C_1$, it follows that both $\phi(v)$ and $\phi(w)$ are within $2\delta$
of the geodesic joining $ax_0$ to $bx_0$.  From this it follows that
$d(\phi(v),\phi(w))$ differs from $d(ax_0,bx_0)-(n+m)C_1$ by at most
$4\delta$.  But since $d(ax_0,bx_0)$ differs by at most $C_1$ from $k
C_1$, we deduce that 
\[(2k-(m+n))C_1 -(C_1+4\delta) \leq d(\phi(v),\phi(w))\leq
(2k-(m+n))C_1 + (C_1+4\delta),\]
from which it immediately follows that
\[C_1 d(v,w) - (2C_1+4\delta)\leq d(\phi(v),\phi(w))\leq C_1 d(v,w)
+(2C_1+4\delta).\] 

In particular, $\phi$ is a $(C_1,2C_1+4\delta)$--quasi-isometric
embedding from the combinatorial horoball $H$ into $X$, and the
theorem is established.
\end{proof}

\section{Rigidity in rank $\geq 2$}\label{s:bigrank}
The purpose of
this section is to establish Theorem \ref{thm:cdr2}, but we will begin
with some lemmas which hold in a slightly broader context.
We suppose that $G$ is a simple Chevalley-Demazure
group scheme of rank at least $2$, that $\Phi$ is a root system for
$G$, and that $R$ is some commutative unital ring containing $\Z$.

\begin{lemma}\label{l:phiprime}
If $\alpha$ and $\beta\in \Phi$,
then there is a $\Phi' = \mathrm{Span}(\Phi')\cap\Phi\subseteq \Phi$
containing $\alpha$ and $\beta$ so that
$\Phi'$ is 
isomorphic to $A_1\times A_1$, $A_2$, $B_2$ or $G_2$.  If $\Phi'\cong
A_1\times A_1$, then $\alpha\neq -\beta$.
\end{lemma}
\begin{proof}
If $\alpha$ and $\beta$ are linearly independent, then
$\Phi'=\mathrm{Span}(\{\alpha,\beta\})\cap\Phi$ is a root system of
rank two, and we simply recall that such a root system is always
isomorphic to one of those listed.

If $\beta=-\alpha$, we assert that there must be some $\gamma\in \Phi$
so that 
$\mathrm{Span}(\{\alpha,\gamma\})\cap\Phi$ is not equal to $A_1\times
A_1$. 
Suppose that there is no such $\gamma$.  Then either 
$\Phi=\overline{\Phi}\times\langle\alpha\rangle$, or
$\Phi=\langle\alpha\rangle$.
Because $G$ is simple $\Phi$ cannot split as a product; because $G$
has rank at least two, $\Phi\neq A_1$.
\end{proof}
For the following two lemmas, we refer to \cite{carter:book}.  Although
the proofs there are done assuming that $R$ is a field, this
assumption is unnecessary;  see also \cite{stein71} or
\cite{steinberg:notes}.
\begin{lemma}\label{l:steinberg}{\em (Steinberg commutator
    relations)}\cite[Theorem 5.2.2]{carter:book}
If $\alpha$, $\beta\in \Phi$ and $t$, $u\in R$, then 
\[[x_\alpha(t),x_\beta(u)]=\prod_{i,j>0\mbox{, }i\alpha+j\beta\in \Phi }
x_{i\alpha+j\beta}(N_{\alpha,\beta,i,j}t_iu_j), \]
where the $N_{\alpha,\beta,i,j}\in \Z$ are integers which depend only
on the order in which the product is taken.
\end{lemma}
\begin{lemma}\label{l:weylconj}\cite[Lemma 7.2.1]{carter:book}
If $\alpha\in \Phi$, $w$ is an element of the Weyl group of $\Phi$,
and $t\in R$,
then 
$x_{w(\alpha)}(t)$ is conjugate either to  $x_\alpha(t)$
or $x_\alpha(-t)$.
\end{lemma}

\begin{lemma}\label{l:distortion}
Let $X$ be a hyperbolic $G(R)$--space.
If $g = x_\alpha(t)$ for some $\alpha\in \Phi$ and $t\in R$, then $g$
does not act hyperbolically on $X$.
\end{lemma}
\begin{proof}
By Lemma \ref{l:phiprime}, there is a subset $\Phi'$ of $\Phi$
containing $\alpha$ which is either isomorphic to $A_2$, $B_2$ or
$G_2$.  In each case, we may apply Lemma \ref{l:steinberg} some number
of times to show that $g$ is distorted in $G$; the details of this are
left to the reader.  By Lemma
\ref{lemma:distortion}, $g$ cannot act hyperbolically on $X$.
\end{proof}

\begin{proposition}\label{p:allfixone}
Let $X$ be a hyperbolic $G(R)$--space,
let $r$, $s\in \Phi$, and let $\rho_1$, $\rho_2\in R$.
Suppose that $p=x_r(\rho_1)$ acts parabolically, fixing some $e\in
\partial X$.  Then $ge=e$ for any other root element $g=x_s(\rho_2)$.
\end{proposition}
\begin{proof}
If $r=s$ or if
$\langle r,s\rangle = A_1\times A_1$, then
$p$ and $g$ commute.  By Lemma \ref{l:commute}, $g(e)=e$, and we are
done.

Otherwise $r$ and $s$ are contained in a two-dimensional root system
$\Phi'\subset \Phi$ which is isomorphic to
$A_2$, $B_2$ or $G_2$, by Lemma \ref{l:phiprime}.

Each case requires a separate argument.  

\begin{case}\label{c:a2}
$\Phi'\cong A_2$.
\end{case}
$A_2 = \{\lambda_i \mid i\in \Z_6\}$ contains six roots, arranged
hexagonally; the angle between $\lambda_i$ and 
$\lambda_j$ is $\frac{|i-j|}{3}\pi$.  Suppose $r=\lambda_i$ and $s=\lambda_j$.
By Lemma \ref{l:weylconj}, $h = x_s(\rho_1)$ is conjugate to either
$p$ or $p^{-1}$, so $h$ is a parabolic, fixing some point
$f\in\partial X$.  
In case $|i-j|=1$, then Lemma \ref{l:steinberg} implies that
$h$ and $p$ commute, and so $f=e$ by Lemma \ref{l:commute}.
In case $|i-j|>1$, one argues by induction on $|i-j|$ to the same
conclusion: $f=e$.

Since $g$ commutes with $h$, we must have $ge=e$, again by Lemma
\ref{l:commute}.

\begin{case}\label{c:b2}
$\Phi'\cong B_2$.
\end{case}

Let $\alpha$ be a short root, and $\beta$ a long root, so that
$\alpha$ and $\beta$ span $\Phi'\subset \Phi$, as in
Figure \ref{f:b_2}.
\begin{figure}[htbp]
\begin{center}
\begin{picture}(0,0)%
\includegraphics{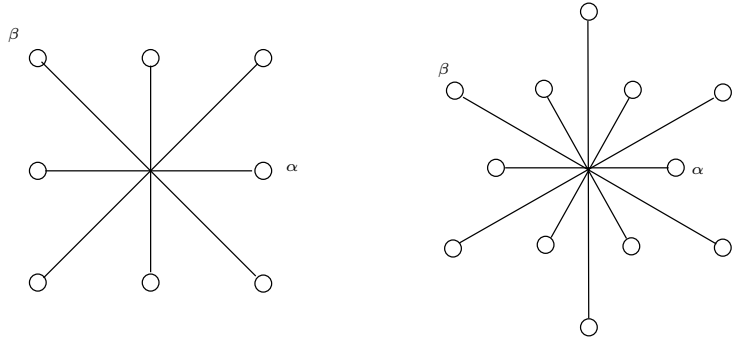}%
\end{picture}%
\setlength{\unitlength}{4144sp}%
\begingroup\makeatletter\ifx\SetFigFont\undefined%
\gdef\SetFigFont#1#2#3#4#5{%
  \reset@font\fontsize{#1}{#2pt}%
  \fontfamily{#3}\fontseries{#4}\fontshape{#5}%
  \selectfont}%
\fi\endgroup%
\begin{picture}(4352,2006)(867,-3021)
\put(4972,-2047){\makebox(0,0)[lb]{\smash{{\SetFigFont{6}{7.2}{\familydefault}{\mddefault}{\updefault}{\color[rgb]{0,0,0}$\alpha$}%
}}}}
\put(3455,-1448){\makebox(0,0)[lb]{\smash{{\SetFigFont{6}{7.2}{\familydefault}{\mddefault}{\updefault}{\color[rgb]{0,0,0}$\beta$}%
}}}}
\put(2546,-2026){\makebox(0,0)[lb]{\smash{{\SetFigFont{6}{7.2}{\familydefault}{\mddefault}{\updefault}{\color[rgb]{0,0,0}$\alpha$}%
}}}}
\put(882,-1238){\makebox(0,0)[lb]{\smash{{\SetFigFont{6}{7.2}{\familydefault}{\mddefault}{\updefault}{\color[rgb]{0,0,0}$\beta$}%
}}}}
\end{picture}%
\caption{$B_2$ and $G_2$.}
\label{f:b_2}
\end{center}
\end{figure}

There are two subcases, depending on whether $r$ is a short or long
root.
\begin{subcase}
The parabolic $p=x_r(\rho_1)$, where $r$ is a short root of $\Phi'$.  
\end{subcase}
Without loss of generality,
we may assume that $r=\alpha$.  If $s\in
\{\alpha,2\alpha+\beta,-\beta\}$, then $g=x_s(\rho_2)$ commutes with
$p$, by Lemma \ref{l:steinberg}. 
Thus by Lemma \ref{l:commute} $ge=e$.

Suppose next that $s = \pm (\alpha+\beta)$.  Then Lemma
\ref{l:steinberg} implies
\begin{equation}\label{c1}
p g p^{-1} g^{-1} = x_{r+s}(N \rho_1\rho_2),
\end{equation}
where $N=N_{r,s,1,1}$ is an integer.  Since $r+s\in
\{2\alpha+\beta,-\beta\}$, we already know that 
$h:=x_{r+s}(N\rho_1\rho_2)$ fixes $e$.  We rearrange \eqref{c1} to give
\begin{equation}
 h^{-1} p = gpg^{-1}.
\end{equation}
Since $p$ and $h$ both fix $e$, so must $gpg^{-1}$.  The element
$gpg^{-1}$ also fixes $g(e)$, since $p$ fixes $e$.
Since $p$ is
parabolic, it can only fix one point in $\partial X$, and so $g(e)=e$.

If $s\in \{\beta, -\alpha\}$, let $r' = \alpha+\beta$; if
$s=-2\alpha-\beta$ let $r' = -\alpha-\beta$.  In any case there is an
element of the Weyl group of $\Phi$ taking $r$ to $r'$; by Lemma
\ref{l:weylconj}, there is a 
$p'= x_{r'}(\pm \rho_1)$ which is
conjugate to $p$ in $E(\Phi,R)$.  Since $p'$ is conjugate to $p$, it
is parabolic; by the argument of the previous paragraph, $p'$ has the
same fixed point as $p$.  If $s\in \{\beta, -2\alpha-\beta\}$, then
$g$ commutes with $p'$, and so $g(e)=e$ by Lemma \ref{l:commute}.
Finally, if
$s=-\alpha$, then we may apply the argument of the previous paragraph
again (with $p'$ and $r'$ in place of $p$ and $r$), to deduce that $g(e)=e$.

\begin{subcase}
The parabolic $p=x_r(\rho_1)$, where $r$ is a long root of $\Phi'$.  
\end{subcase}
In this case, we may assume for instance that $r=2\alpha+\beta$.  If
$s\in\{2\alpha+\beta, \alpha+\beta, \alpha, \beta, -\beta\}$, then
Lemma \ref{l:steinberg} implies that $g=x_s(\rho_2)$ commutes with $p$,
and so $g(e)=e$ by Lemma \ref{l:commute}.

Suppose then that $s\in \{-\alpha,-2\alpha-\beta,-\alpha-\beta\}$.
If $s\in\{-\alpha,-2\alpha-\beta\}$, let $r'=\beta$; if $s=
-\alpha-\beta$, then let $r'=-\beta$.  In either case, there is an
element of the Weyl group taking $r$ to $r'$, and so there is a
parabolic $p' = x_{r'}(\pm \rho_1)$ conjugate to $p$ by Lemma
\ref{l:weylconj}.   By the previous paragraph, $p'$ has the same fixed
point as $p$ does.  Applying the previous paragraph with $p'$ and $r'$
in the place of $p$ and $r$ implies that $g(e)=e$ for $g =
x_s(\rho_2)$. This completes the proof of Case \ref{c:b2}.

\begin{case}\label{c:g2}
$\Phi'\cong G_2$.
\end{case}
Let $\alpha$ be a short root, and $\beta$ a long root, so that
$\alpha$ and $\beta$ span $\Phi'\subset \Phi$, as in
the right half of Figure \ref{f:b_2}.

Again there are two subcases, depending on whether $r$ is a short or long
root.

\begin{subcase}
The parabolic $p=x_r(\rho_1)$, where $r$ is a long root of $\Phi'$.  
\end{subcase}
If $s\in \{-\beta, 3\alpha+\beta, \pm(3\alpha+2\beta)\}$, then Lemma
\ref{l:steinberg} implies that $g=x_s(\rho_2)$ commutes with $p$, and
so $g(e)=e$ by Lemma \ref{l:commute}.

Suppose that $s\in \{2\alpha+\beta,-\alpha-\beta\}$.  Lemma
\ref{l:steinberg} implies that 
\begin{equation}
p g p^{-1} g^{-1} = x_{r+s}(N\rho_1\rho_2)=: h,
\end{equation}
for some integer $N$.  Exactly as in Case \ref{c:b2}, $h$ commutes
with $p$, and so $h(e)=e$.  Thus $e= h^{-1} p (e) = gpg^{-1}(e)$ and
the parabolic $gpg^{-1}$ fixes $e$.  Again since $gpg^{-1}$ also fixes
$g(e)$, we must have $g(e)=e$.

Using Lemma \ref{l:weylconj} repeatedly we discover that for every
short root $r'$ there is a parabolic element $p'=x_{r'}(\pm \rho_1)$
with $p'(e)=e$.  Since $g$ must commute with some such element,
$g(e)=e$ as well.
\begin{subcase}
The parabolic $p=x_r(\rho_1)$, where $r$ is a long root of $\Phi'$.  
\end{subcase}
Without loss of generality, we may assume that $r = 3\alpha+2\beta$.

If the inner product of $s$ with $r$
is nonnegative, then $s\in\{ \pm \alpha, \beta, \alpha+\beta,
\alpha+2\beta,\alpha+3\beta, 3\alpha+2\beta\}$, and $g$ commutes with
$p$ by Lemma \ref{l:steinberg}, and so $g(e)=e$.  

Otherwise, a (possibly
repeated) application of Lemma \ref{l:weylconj} implies that $g$
commutes with a parabolic $p' = x_{r'}(\pm \rho_1)$ for some long root
$r$ of $\Phi'$, and with $p'(e)=e$.  This again implies via Lemma
\ref{l:commute} that $g(e)=e$.  This completes the proof in Case \ref{c:g2}.
\end{proof} 

We are now ready to give the proof of Theorem \ref{thm:cdr2}.  
\begin{proof}
A result of Tavgen$'$ \cite{tavgen:bg} shows that $G(\mathcal{O})$ is
boundedly generated by its root subgroups.
The ring of integers $\mathcal{O}$ is finitely generated as an Abelian
group; choose generators $\mu_1,\ldots,\mu_k$.  It follows from
Tavgen$'$'s result that the set
\[ S=\{x_\alpha(\mu_i) \mid \alpha \in \Phi\mbox{, }1\leq i\leq k\}\]
boundedly generates $G(\mathcal{O})$.

Each of these generators acts hyperbolically, elliptically or
parabolically on $X$.  
By Lemma \ref{l:distortion}, none can act hyperbolically.
If all of the root elements act elliptically, then it follows from bounded
generation that the orbit of a point under the action of $G$ must be
bounded.

We therefore may assume that some $x_\alpha(\mu_i)$ acts parabolically on
$X$, fixing a single point $e\in \partial X$.  It follows from
Proposition \ref{p:allfixone} that all the root subgroups will fix
this point $e$, and so $G(\mathcal{O})$ fixes $e$.

By Proposition \ref{p:pseudo}, the pseudocharacter 
$p_\mathbf{x}\co G(\mathcal{O})\to \R$ determined by a sequence
$\mathbf{x}=\{x_i\}$ tending to $e$
is nonzero exactly on the hyperbolic elements.  Thus $p_\mathbf{x}(g)=0$
whenever $g$ lies in a root subgroup.  

An elementary argument (see for example \cite[Proposition
  5]{kotschick:quasi}) shows that a pseudocharacter $p$ on a boundedly
generated group 
is determined by its values on
the bounded generators;  thus $p\equiv 0$ on $G(\mathcal{O})$.  It
follows that no element of $G(\mathcal{O})$ acts hyperbolically on
$X$. 
By Theorem
\ref{t:horeq}, the $G(\mathcal{O})$--space $X$ has an invariant horoball.
\end{proof}

\section{Remarks on rank one}\label{section:rankone}
One can also ask what actions rank one Chevalley groups have on
hyperbolic spaces.  If $\mc{O}$ is a number ring with finitely many
units, then $SL(2,\mc{O})$ is a lattice either in $SL(2,\R)$ or
$SL(2,\C)$.  In particular, it has a proper non-elementary action on
$\H^2$ or $\H^3$.  Moreover, such a group admits uncountably many
distinct pseudocharacters (AKA homogeneous quasi-(homo)morphisms) up
to scale \cite{fujiwara:gromovhyperbolic,bestvinafujiwara:mcg}.  Each
such ``projective pseudocharacter'' gives rise to a 
quasi-action on $\R$; no two such are equivalent.  Moreover, these
often give rise to quasi-actions on more complicated trees
\cite{manning:cocycles}.
The groups $SL(2,\mc{O})$ where $\mc{O}$ has infinitely many units are
more rigid.
In this section we
apply the main result of \cite{manning:qfa} to the special case of
actions on quasi-trees (defined below), and
speculate on the general situation.

Recall that a group $G$ is said to have 
Property \FA\ if every action by $G$ on a simplicial tree $T$ has a
fixed point.
\begin{definition}
A \emph{quasi-tree} is a graph which is quasi-isometric to a tree.
\end{definition}
\begin{definition}
A group $G$ has property \QFA\ if every action by $G$ on a quasi-tree
$X$ has a bounded orbit. 
\end{definition}
\begin{remark}
This is differently worded than the definition in \cite{manning:qfa},
but easily seen to be equivalent.  Note that \QFA\ implies \FA, but
not vice versa.
\end{remark}
As quasi-trees are in particular Gromov hyperbolic spaces which admit
no parabolic isometries (see Section 3.2 of \cite{manning:qfa}),
Theorem \ref{thm:cdr2} implies that higher rank Chevalley
groups over number rings have property \QFA.

We recall a definition and a theorem from \cite[Section 4]{manning:qfa}.

\begin{definition}
Let $G$ be a group, and let $g$ be an element of $G$.  We will say $g$
is a \emph{stubborn element of $G$} if for all
$H< G$ with $[G:H]\leq 2$, there exists some integer $k_H>0$ so that
$g^{k_H}\in [H,H]$.  
\end{definition}

\begin{theorem}\label{thm:qfa}\cite[Theorem 4.4]{manning:qfa}
Let $G$ be a group which is boundedly generated by elements 
$g_1,\ldots,g_n$, so
that for each $i$, $g_i$ is a stubborn element of $B_i$ for some
amenable $B_i < G$. Then $G$ has Property \QFA.
\end{theorem}
Note that the above theorem was misstated slightly in
\cite{manning:qfa}; the
word ``amenable'' was inadvertently omitted.

Here's an easy lemma:
\begin{lemma}\label{lemma:finiteindex}
Let $H<G$ be a subgroup of finite index.  If $H$ has Property
\QFA, then so does $G$.
\end{lemma}

Our methods in the higher rank case use heavily the bounded generation of
$G(\mathcal{O})$ established by Tavgen$'$ in \cite{tavgen:bg}
for Chevalley groups over rings of integers of algebraic
number fields.  There is an analogous result of Carter, Keller, and Paige
in rank $1$, at least for $SL(2,\cdot)$ and certain number
fields:
\begin{theorem}\label{thm:ckp2}
\cite{wittemorris:ckp}  
For any integer $d>1$ there is an $r=r(d)$ so that the following is true.
Let $K$ be a number field of degree $d$ over $\Q$, and let
$\mathcal{O}$ be the ring of integers of $K$.  If the $\mc{O}$ has
infinitely many units, then:
\begin{enumerate}
\item every element of $E(2,\mathcal{O})$ is a product of at most $r$
  elementary matrices, and
\item the index of $E(2,\mathcal{O})$ in $SL(2,\mathcal{O})$ is at
  most $r$.  
\end{enumerate}
\end{theorem}
In the above statement, $E(2,\mathcal{O})$ is the subgroup of
$SL(2,\mathcal{O})$ generated by the root subgroups (the strictly
upper triangular and strictly lower triangular matrices).
The following statement implies Theorem \ref{thm:cdr1}:
\begin{theorem}\label{thm:rankone}
If $\mathcal{O}$ is the ring of integers of an algebraic number field
and $\mathcal{O}$ has infinitely many units, then $SL(2,\mathcal{O})$
has Property \QFA.
\end{theorem}
\begin{proof}
By Lemma \ref{lemma:finiteindex} it suffices to show that
$E(2,\mathcal{O})$ has property \QFA.
If $\Lambda= \{\lambda_1,\ldots,\lambda_n\}$ is an integral basis
for the number field $\mathcal{O}$, then
Theorem \ref{thm:ckp2} implies that $E(2,\mathcal{O})$ is boundedly
generated by the $2n$ elements $\textMtwo{1}{\lambda_i}{0}{1}$ and
$\textMtwo{1}{0}{\lambda_i}{1}$.  
\begin{claim}
For each $i\in\{1,\ldots,n\}$,
$\textMtwo{1}{\lambda_i}{0}{1}$ is a stubborn element of 
$B=\textMtwo{*}{*}{0}{*}\cap E(2,\mathcal{O})$.
\end{claim}
\begin{proof}
On page 189 of \cite{carter:book}, Carter observes that
$\textMtwo{t}{0}{0}{t^{-1}}$ can be written as a product of elementary
matrices, for any invertible $t\in \mc{O}$, as follows.  
If
$\lambda$ is any invertible element, we may write
\[\Mtwo{0}{\lambda}{-\lambda^{-1}}{0} =
\Mtwo{1}{\lambda}{0}{1}\Mtwo{1}{0}{-\lambda^{-1}}{1}\Mtwo{1}{\lambda}{0}{1},\]
and then note that
\begin{equation}\label{computation}
\Mtwo{t}{0}{0}{t^{-1}}=\Mtwo{0}{t}{-t^{-1}}{0}\Mtwo{0}{-1}{1}{0}.
\end{equation}
Carter's observation shows that $\textMtwo{\omega}{0}{0}{\omega^{-1}}$ is in
$E(2,\mc{O})$ for any unit $\omega$ in $\mc{O}$.  
Computing the commutator of \textMtwo{\omega}{0}{0}{\omega^{-1}} and
\textMtwo{1}{\lambda}{0}{1} for $\omega$ a unit of $R$ and $\lambda\in R$
yields:
\begin{equation}\label{computation2}
\left[\Mtwo{\omega}{0}{0}{\omega^{-1}},\Mtwo{1}{\lambda}{0}{1}\right]
 = \Mtwo{1}{(1-\omega^2)\lambda}{0}{1}.
\end{equation}

By assumption, the group of units of $\mc{O}$ is infinite.
Dirichlet's units theorem (see, e.g. \cite[Appendix B]{stewarttall})
implies that we may
choose $\omega_0\in \mc{O}^*$ a unit of infinite order. 
Let $H<B$ be a subgroup of index at most two.  Then $H$ must contain
$\textMtwo{\omega_0^2}{0}{0}{\omega_0^{-2}}$ and
$\textMtwo{1}{2 r}{0}{1}$ for all $r\in \mc{O}$.  It follows from the
computations \eqref{computation} and \eqref{computation2}
that $\textMtwo{1}{i}{0}{1}\in [H,H]$ for all
$i$ in the ideal $I$ generated by $2(1-\omega_0^2)$.  Let $N$ be the
order of $R/I$. (The number $N$ is also called the \emph{norm} of $I$;
that it is finite when $I\neq (0)$ is an elementary fact of algebraic
number theory;  see, e.g., \cite[Chapter 5]{stewarttall}.)
For any of the $\lambda_i$, we have
$\textMtwo{1}{\lambda_i}{0}{1}^N=\textMtwo{1}{N\lambda_i}{0}{1}\in [H,H]$, and
so $\textMtwo{1}{\lambda_i}{0}{1}$ is stubborn.
\end{proof}
It remains to observe that $B<E(2,\mathcal{O})$ is solvable, and hence
amenable.
We may now apply Theorem \ref{thm:qfa} to conclude that
$E(2,\mathcal{O})$ has Property \QFA.
\end{proof}
\begin{remark}
It was already known \cite[p. 68]{serre:trees} that the groups covered
by Theorem \ref{thm:rankone} possessed property \FA.
\end{remark}
\begin{remark}
Another proof of \ref{thm:rankone} may be given as follows:  First
show that every unipotent is distorted.  It follows that the bounded
generators cannot (quasi)-act hyperbolically.  It is shown in
\cite[Corollary 3.6]{manning:qfa} that there are no parabolic isometries of
quasi-trees, and so each of the bounded generators (quasi)-acts
elliptically.  It then follows from bounded generation that any orbit
is bounded.
\end{remark}

Finally, we speculate on the variety of hyperbolic $\Gamma$-spaces, for
$\Gamma=SL(2,\mathcal{O})$, where $\mathcal{O}$ is the
ring of integers of a number field $k$.  
We have already remarked that $\Gamma$ is a lattice in 
\[\prod_{i=1}^s SL(2,\R)\times \prod_{i=1}^t SL(2,\C),\]
where $s$ and $t$ are the number of real and complex places
respectively.  Projection to some factor gives an isometric action
either on $\H^2$ or $\H^3$.  Call a hyperbolic $\Gamma$--space
\emph{standard} if it is equivalent to $\H^2$ or $\H^3$ with one of
these actions.

\begin{conjecture}\label{conjecture:realquadratic}
Every quasi-action by $\Gamma$ on a Gromov hyperbolic metric space
either has an invariant horoball or is standard.
\end{conjecture}

%
%
%
%

\def\cprime{$'$} \providecommand\url[1]{\texttt{#1}}

\end{document}

%% file: bar.pstex_t
\begin{picture}(0,0)%
\includegraphics{bar.pstex}%
\end{picture}%
\setlength{\unitlength}{3947sp}%
\begingroup\makeatletter\ifx\SetFigFont\undefined%
\gdef\SetFigFont#1#2#3#4#5{%
  \reset@font\fontsize{#1}{#2pt}%
  \fontfamily{#3}\fontseries{#4}\fontshape{#5}%
  \selectfont}%
\fi\endgroup%
\begin{picture}(4707,2311)(1082,-1752)
\put(1108,425){\makebox(0,0)[lb]{\smash{{\SetFigFont{9}{10.8}{\familydefault}{\mddefault}{\updefault}{\color[rgb]{0,0,0}$\uu{v}_1$}%
}}}}
\put(2367,425){\makebox(0,0)[lb]{\smash{{\SetFigFont{9}{10.8}{\familydefault}{\mddefault}{\updefault}{\color[rgb]{0,0,0}$\uu{w}_1$}%
}}}}
\put(1165,-1693){\makebox(0,0)[lb]{\smash{{\SetFigFont{9}{10.8}{\familydefault}{\mddefault}{\updefault}{\color[rgb]{0,0,0}$\uu{v}_2$}%
}}}}
\put(2482,-1693){\makebox(0,0)[lb]{\smash{{\SetFigFont{9}{10.8}{\familydefault}{\mddefault}{\updefault}{\color[rgb]{0,0,0}$\uu{w}_2$}%
}}}}
\put(4276,425){\makebox(0,0)[lb]{\smash{{\SetFigFont{9}{10.8}{\familydefault}{\mddefault}{\updefault}{\color[rgb]{0,0,0}$\overline{\uu{v}_1}$}%
}}}}
\put(5530,425){\makebox(0,0)[lb]{\smash{{\SetFigFont{9}{10.8}{\familydefault}{\mddefault}{\updefault}{\color[rgb]{0,0,0}$\overline{\uu{w}_1}$}%
}}}}
\put(4333,-1683){\makebox(0,0)[lb]{\smash{{\SetFigFont{9}{10.8}{\familydefault}{\mddefault}{\updefault}{\color[rgb]{0,0,0}$\overline{\uu{v}_2}$}%
}}}}
\put(5644,-1683){\makebox(0,0)[lb]{\smash{{\SetFigFont{9}{10.8}{\familydefault}{\mddefault}{\updefault}{\color[rgb]{0,0,0}$\overline{\uu{w}_2}$}%
}}}}
\end{picture}%

%% file: deep.pstex_t
\begin{picture}(0,0)%
\includegraphics{deep.pstex}%
\end{picture}%
\setlength{\unitlength}{3947sp}%
\begingroup\makeatletter\ifx\SetFigFont\undefined%
\gdef\SetFigFont#1#2#3#4#5{%
  \reset@font\fontsize{#1}{#2pt}%
  \fontfamily{#3}\fontseries{#4}\fontshape{#5}%
  \selectfont}%
\fi\endgroup%
\begin{picture}(1171,2766)(961,-2305)
\put(1276,314){\makebox(0,0)[lb]{\smash{{\SetFigFont{12}{14.4}{\familydefault}{\mddefault}{\updefault}{\color[rgb]{0,0,0}$a$}%
}}}}
\put(2101,314){\makebox(0,0)[lb]{\smash{{\SetFigFont{12}{14.4}{\familydefault}{\mddefault}{\updefault}{\color[rgb]{0,0,0}$b$}%
}}}}
\put(1876,-1861){\makebox(0,0)[lb]{\smash{{\SetFigFont{12}{14.4}{\familydefault}{\mddefault}{\updefault}{\color[rgb]{0,0,0}$x_{i(b,m)}$}%
}}}}
\put(1726,-2236){\makebox(0,0)[lb]{\smash{{\SetFigFont{12}{14.4}{\familydefault}{\mddefault}{\updefault}{\color[rgb]{0,0,0}$x_I$}%
}}}}
\put(976,-1336){\makebox(0,0)[lb]{\smash{{\SetFigFont{12}{14.4}{\familydefault}{\mddefault}{\updefault}{\color[rgb]{0,0,0}$x_{i(a,n)}$}%
}}}}
\put(1951,-736){\makebox(0,0)[lb]{\smash{{\SetFigFont{12}{14.4}{\familydefault}{\mddefault}{\updefault}{\color[rgb]{0,0,0}$\phi(w)$}%
}}}}
\put(976,-436){\makebox(0,0)[lb]{\smash{{\SetFigFont{12}{14.4}{\familydefault}{\mddefault}{\updefault}{\color[rgb]{0,0,0}$\phi(v)$}%
}}}}
\end{picture}%